\documentclass[11pt]{article}
\usepackage[utf8]{inputenc}
\usepackage[T1]{fontenc}
\usepackage{amsmath,amssymb,amsthm,bbm,float,graphicx,geometry,lmodern,mathtools,parskip,setspace,subcaption}
\usepackage[colorlinks, citecolor = red!80!black, linkcolor = blue!80!black, breaklinks, pdfauthor ={"Peter Gracar, Lukas L\"{u}chtrath, Christian M\"{o}nch"}]{hyperref}
\usepackage{mathrsfs}
\usepackage{enumerate}
\usepackage{nicefrac}
\usepackage{todonotes}
\usepackage{nth}
\usepackage{titling}
\usepackage{orcidlink}

\usepackage{tikz}
\usetikzlibrary{backgrounds}
\usetikzlibrary{patterns}
\usetikzlibrary{positioning, shapes.geometric}
\usetikzlibrary{external, arrows, decorations.pathmorphing} 

\usepackage{caption}
\usepackage{graphicx}
\usepackage{todonotes}
\graphicspath{{Images/}}
\allowdisplaybreaks
\usepackage[raggedright]{titlesec} 
\usepackage{booktabs} 
\DeclarePairedDelimiter{\set}{\lbrace}{\rbrace}
\newcommand{\N}{\mathbb{N}}

\newtheorem{theorem}{Theorem}[section]
\newtheorem{lemma}[theorem]{Lemma}
\newtheorem{corollary}[theorem]{Corollary}
\newtheorem{proposition}[theorem]{Proposition}

\theoremstyle{definition}
\newtheorem{definition}[theorem]{Definition}
\newtheorem{remark}[theorem]{Remark}

\renewcommand{\P}{\mathbb{P}}
\newcommand{\x}{\mathbf{x}}
\newcommand{\y}{\mathbf{y}}
\newcommand{\1}{\mathbbm{1}}
\newcommand{\G}{\mathscr{G}}
\newcommand{\C}{\mathscr{C}}
\newcommand{\X}{\mathbf{X}}
\newcommand{\E}{\mathbb{E}}
\newcommand{\0}{\mathbf{0}}
\newcommand{\R}{\mathbb{R}}
\renewcommand{\d}{\mathrm{d}}
\newcommand{\z}{\mathbf{z}}

\newcommand{\Z}{\mathbb{Z}}

\newcommand{\Vl}{\mathscr{V}_{\ell}^n(\vartheta^*)}
\newcommand{\Vr}{\mathscr{V}_{r}^n(\vartheta^*)}
\newcommand{\cG}{\G}
\newcommand{\cL}{\mathcal{L}}
\newcommand{\cR}{\mathcal{R}}
\newcommand{\cE}{\mathcal{E}}
\newcommand{\cF}{\mathcal{F}}
\newcommand{\cA}{\mathcal{A}}
\newcommand{\ER}{\operatorname{ER}}

\usepackage[style = numeric, sorting=nyt, url = false, abbreviate=false, maxbibnames=9, sortcites=true, doi = true, backend = biber, giveninits = true, isbn = false]{biblatex}
\renewbibmacro{in:}{\ifentrytype{article}{}{\printtext{\bibstring{in}\intitlepunct}}}
\DeclareFieldFormat{journaltitle}{\mkbibemph{#1\isdot}}
\bibliography{infiniteComponents.bib}

\pretitle{\centering\LARGE\scshape}
 \posttitle{\vskip 0.75cm}

 \predate{\vskip 0.75 cm \centering\large}
 \postdate{\par}

\thanksmarkseries{arabic}

\author{Peter Gracar \orcidlink{0000-0001-8340-8340} \thanks{University of Leeds, Woodhouse, Leeds LS2 9JT, United Kingdom} \\P.Gracar@leeds.ac.uk
	\and
	Lukas L\"{u}chtrath \orcidlink{0000-0003-4969-806X} \thanks{Weierstrass Institute for Applied Analysis and Stochastics, Mohrenstraße 39, 10117 Berlin, Germany} \\ lukas.luechtrath@wias-berlin.de
	\and
	Christian M\"onch \orcidlink{0000-0002-6531-6482} \thanks{Johannes Gutenberg University Mainz, Saarstraße 21, 55122 Mainz, Germany} \\
	cmoench@uni-mainz.de
}

\title{Finiteness of the percolation threshold for
		inhomogeneous long-range models in
		one dimension} 
        
\date{August 15, 2025}

\usepackage{etoolbox}
\ifundef{\abstract}{}{\patchcmd{\abstract}%
{\quotation}{\quotation\noindent\ignorespaces}{}{}}
\begin{document}
\maketitle 
\begin{spacing}{0.9}
\begin{abstract} 
We consider inhomogeneous spatial random graphs on the real line. Each vertex carries an i.i.d.\ weight and edges are drawn such that short edges and edges to vertices with large weights occur with higher probability. This allows the study of models with long-range effects and heavy-tailed degree distributions. We introduce a new coefficient \(\delta_\textup{eff}\) which quantifies the influence of heavy-tailed degrees on long-range connections. We show that \(\delta_\textup{eff}<2\) is sufficient for the existence of a supercritical percolation phase in the model and that \(\delta_\textup{eff}>2\) always implies the absence of percolation. In particular, our results complement those in Gracar et al. (Adv.\ Appl.\ Prob., 2021), where sufficient conditions were given for the soft Boolean model and the age-dependent random connection model for both the existence and the absence of a subcritical percolation phase. Our results further provide a criterion for the existence or non-existence of a giant component in large finite graphs. 

\medskip\noindent
{\normalsize\textbf{AMS-MSC 2020}:} 05C80 (Primary), 60K35 (Secondary)
\end{abstract}
\end{spacing}
\begin{spacing}{1}

\newpage


\section{Introduction and overview of the results}\label{SecIntro}

\subsection*{Background and motivation}
We study a large class of spatially inhomogeneous Bernoulli-type percolation models collectively known as the \emph{weight-dependent random connection model} that was first described in \cite{GHMM2022} and further analysed in \cite{GLM2021,GGM21}.
The model is an extension of classical spatial percolation models in the same sense as the non-spatial inhomogeneous random graph model introduced by Bollob\'as et al.\ in \cite{BollobasJansonRiordan2007} generalises the Erd\H{o}s-R\'enyi random graph. Very broadly speaking, the framework we work in is as follows: Consider a family of graphs \((\G_\beta\colon\beta\geq 0)\) on a countably infinite vertex set embedded into \(\mathbb{R}^d\). Here, \(\beta\) is some intensity parameter that controls connectivity in the sense that increasing \(\beta\) results in more connections on average and $\beta=0$ corresponds to a graph without any edges. Common assumptions in this context are 
\begin{enumerate}[(A)]
	\item given the vertex set, edges occur independently of one another whenever they do not share a common vertex;
	\item the vertex set and the connection probabilities are translation invariant\footnote{Since we cover both lattice models and models based on ergodic point processes we do not specify precisely the group of shifts involved, this will always be clear from the context.};
	\item spatially close vertices are more likely to be connected by an edge than far apart vertices.
\end{enumerate}
Note that assumption (A) covers i.i.d.\ Bernoulli site percolation as well as i.i.d.\ bond percolation on lattices. Typically, the vertex set in the models we have in mind is either given by \(\Z^d\) or by a homogeneous Poisson process on $\R^d$, the case we focus on in this section. However, let us mention that our results hold for a large class of renewal processes, see Definition~\ref{DefEvenlySpaced} and Remark~\ref{rem:evenly}. The classical subject of percolation theory is the following emergence phenomenon: Is there a critical intensity parameter \(\beta_c\in(0,\infty)\) such that the graph \(\G_\beta\) contains an infinite connected component, or \emph{infinite cluster}, for \(\beta>\beta_c\) and no infinite cluster for \(\beta<\beta_c\)? Assumptions (A) and (B) generally entail that the existence of an infinite cluster is a $0$-$1$ event. Moreover, if an infinite cluster exists, it is almost surely unique under very mild assumptions on the distribution of $\G_\beta$ \cite{BurtonKeane89,JacobMoerters2017,chebunin2024uniqueness}.

This article is almost exclusively devoted to dimension $d=1$ and \(\G_\beta\) with finite mean degree. In the one-dimensional case, both in classical long-range Bernoulli percolation models and in Poisson-Boolean percolation (which are both special instances of the weight-dependent random connection model) we have that $\beta_c=0$ if and only if \(\G_\beta\) has infinite mean degree \cite{NewmanSchulman1986,MeesterRoy1996}. This is not the case for the weight-dependent random connection model in general and in recent years several papers have identified parameter regimes for which $\beta_c=0$ whilst \(\G_\beta\) has finite mean degree for various instances of the model with heavy-tailed degree distributions, e.g.\ \cite{Yukich2006, DeijfenHofstadHooghiemstra2013,DeprezWuthrich2019,GLM2021}.

In this work we address the complementary question whether $\beta_c<\infty$. The topological restrictions of the line make it rather hard to generate infinite clusters and many `natural' models admit no percolation phase. In nearest-neighbour percolation on \(\Z\), either all edges/vertices are present or the graph contains only finite components~\cite{Grimmett1999}. Similarly, in the one dimensional Gilbert graph~\cite{Gilbert61}, all components are finite as well~\cite{MeesterRoy1996}. Generally, if each vertex can only have neighbours within bounded distance and if there is a modicum of independence in the connection probabilities then there will be no percolation. Moreover, the absence of a percolation phase extends to some models with unbounded interaction range, such as the Poisson-Boolean model with integrable radius distribution, 
even though the degree distribution of \(\G_\beta\) might be heavy-tailed in some instances, see \cite{Hall85, MeesterRoy1996, Gouere08, JahnelTobiasCali2022}. Conversely, if $d\ge 2$ then \(\G_\beta\) tends to contain an infinite component for large \(\beta\) in many cases, which is essentially a consequence of the existence of a supercritical percolation phase in nearest-neighbour Bernoulli percolation on $\Z^2$~\cite{Grimmett1999,MeesterRoy1996}.

The picture for \emph{long-range Bernoulli bond percolation} on $\Z$ is different. Here, any pair of vertices \(x\) and \(y\) is connected independently with probability proportional to 
\[
1\wedge \beta |x-y|^{-\delta} \quad \text{ for some }\delta>0,
\]
we denote this model by $\operatorname{LRP}(\delta)$. The case $\delta\leq1$ yields graphs that have infinite expected degree and shall not interest us here. If, however, $\delta\in(1,2)$ then we have \(0<\beta_c<\infty\) and if \(\delta>2\), then no infinite cluster exists for any $\beta>0$ \cite{Schulman1983,NewmanSchulman1986,DCGT2019}. Our results suggest that the behaviour of $\operatorname{LRP}(\delta)$ in these two regimes is paradigmatic for those instances of the weight-dependent random connection model that show genuine long-range interactions. The third regime, which for classical long-range percolation on $\Z$ coincides with the boundary case $\delta=2$ is more delicate. In fact, it is known that \(\beta_c<\infty\) in $\operatorname{LRP}(2)$ and that there even exists an infinite cluster at the critical point $\beta=\beta_c$ \cite{AizenmanNewman1986,DCGT2019}, which is rather atypical. We refer to this particular model as \emph{scale-invariant long-range percolation}.

Our arguments rely on a renormalisation scheme, developed in Section~\ref{secInfiniteComponent}, which is tailored to graphs with positively correlated edges and with vertex locations given by a renewal process. A similar approach was used by Duminil-Copin et al.\ in~\cite{DCGT2019} to study \(\operatorname{LRP}(2)\) on $\Z$. The positive correlations between edges are built into our model by way of adding weights to the vertices. Dealing with the arising correlations is one of the key mathematical contributions of this paper. There are two major issues that have to be overcome:  First, we have to control the effect of the vertex weights on the number of connections of far apart sets of vertices, cf.\ Lemma~\ref{lemLongConnection}. Secondly, we have to deal with the existence of well-connected clusters when going through the scales of the renormalisation. The information of existence of such clusters provides information about the connectivity in the model and thus influences the distributions of the vertex weights and edge probabilities, a problem which always arises when employing multi-scale arguments to inhomogeneous percolation models. In Lemma~\ref{lem:Decorrelate}, we show how to decorrelate the corresponding probabilities involved, using the FKG-inequality on suitably chosen sub-$\sigma$-fields. Our method formalises the intuition that a vertex in a large cluster should not have fewer neighbours than a `typical' vertex and was used in~\cite{Monch2023} to obtain a finite size criterion for percolation in inhomogeneous long-range percolation models.

We apply our main result to a large variety of model instances from the literature, cf.\ Proposition~\ref{lem:InterPol}. Combined with results of~\cite{GLM2021,DeprezWuthrich2019}, this allows us to establish an almost complete phase diagram summarising the existence of percolation thresholds in several model instances of interest, see Figure~\ref{figPercolation} and Theorem~\ref{propPercolation}. Finally, our main result yields information about the occurrence of giant components in finite version of our model, cf.~Corollary~\ref{cor:finiteGraphs}.   

Throughout the paper we use the following notation for non-negative functions: we write \(f=o(g)\) as \(x\to\infty\) if \(\lim_{x\to\infty} \nicefrac{f(x)}{g(x)}=0\), and \(f\asymp g\) if \(f/g\) is asymptotically bounded from zero and infinity. The number of elements in a finite set \(A\) is denoted as \(\sharp A\). We further indicate by \(a-\) and by \(a+\) directional limits towards $a\in\R$. 

\subsection{Discussion of results}
Our model has two principal components:
\begin{itemize}
	\item the \emph{kernel}, a symmetric function \(g\colon(0,1)^2\to(0,\infty)\) which is non-decreasing in both arguments and satisfies $$\int_0^1\int_0^1 \frac{1}{g(s,t)}\,\d s\,\d t<\infty;$$
	\item the \emph{profile}, a non-increasing and integrable function \(\rho\colon(0,\infty)\to[0,1]\).
\end{itemize}

For given $\beta>0$, we now generate \(\G_\beta\) in two steps. First, we sample vertices. In this introductory section, the \emph{vertex locations} are always distributed according to a unit intensity Poisson process $\eta$ on $\R$, we discuss more general choices later on in Section~\ref{sec:formalConstr}. Each vertex location $x\in\eta$ is endowed with an independent $\operatorname{Uniform(0,1)}$ \emph{vertex mark} $t_x$. The pairs $\x=(x,t_x)$ form the vertices of $\G_\beta$ and we use the notation $\mathscr{X}$ for the collection of all vertices. Secondly, given $\mathscr{X}$, we connect any two vertices $\x,\y\in\mathscr{X}$ independently by an edge with probability
\[
\rho\Big(\tfrac{1}{\beta}\,{g(t_x,t_y)\,|x-y|}\Big),
\]
and we denote the presence of an edge between $\x$ and $\y$ by $\x\sim\y.$ Due to our monotonicity assumptions on $\rho$ and $g$, short edges and connections to vertices with small marks are more probable than long edges and edges between vertices of high mark. One can think of the marks as inverse vertex weights, which give the weight-dependent random connection model its name. The integrability assumptions on $\rho$ and $g$ assure that $\G_\beta$ has finite mean degree. 
The parameter \(\beta>0\) directly affects the edge probability but can also be interpreted as controlling the density of Poisson points, since it rescales distances. In fact, our approach enables us to work with more general vertex locations than a Poisson process and we give a detailed construction of the model in Section~\ref{SecGen}.
By construction, the kernel \(g\) determines the degree distribution of the model up to a multiplicative constant, whereas the profile specifies the geometric restrictions of the model \cite{GHMM2022, GGLM2019}. 

In the literature, two types of profile functions have attracted the most attention: profile functions with bounded support, e.g.\ the indicator \(\rho=\mathbbm{1}_{[0,1]}\), and profile functions with polynomial decay of index \(\delta>1\), i.e., if there exist constants \(0<c<C\) such that
\begin{equation}\label{eq:rhodecay}
	c (1\wedge r^{-\delta})\leq \rho(r)\leq C(1\wedge r^{-\delta}), \ \forall r>0.
\end{equation}
Alternatively, one could also assume \(\rho\) to be regularly varying with index \(-\delta\). An application of standard Potter bounds \cite{Bingham1989} shows that there is no qualitative change in the behaviour of existence of infinite clusters compared to the choice~\eqref{eq:rhodecay} in any of the concrete examples discussed in this paper. The decay index \(\delta\) measures how strict the geometric restrictions of the model are. Decreasing values of \(\delta\) lead to more and more long edges. On the contrary, the larger the value of \(\delta\), the stricter the geometric restrictions are, the extreme case being a profile \(\rho\) of bounded support. We therefore refer to models constructed with a profile of bounded support as \emph{hard} models and to models with a polynomial decay of index \(\delta\in(1,\infty)\) as \emph{soft} models. Profiles of unbounded support but with lighter tails than polynomial are included in our approach but show no qualitatively different behaviour in the existence of infinite components compared to the bounded support case, cf.~ \cite{GGM21, Steif1996}. This article is primarily devoted to soft models with heavy-tailed degree distributions but our main results hold for any choices of \(\rho\) and \(g\) satisfying the monotonicity and integrability conditions stated above.

By tuning $\rho$ and $g$ one obtains (versions of) many previously studied percolation models such as the previously mentioned Poisson-Boolean model, the random connection model \cite{MeesterRoy1996}, scale-free percolation \cite{DeijfenHofstadHooghiemstra2013,DeprezWuthrich2019} and the age-dependent random connection model \cite{GGLM2019}. A further closely related model class of similar generality is that of geometric inhomogeneous random graphs \cite{BringmannKeuschLengler2019,HofstadHoornMaitra2021, JorritsmaEtAl2023}. A more comprehensive discussion of related models is given in \cite{GHMM2022}.

\subsubsection*{Existence of percolation threshold via effective decay exponents.} We now formulate our main result which relates finiteness of $\beta_c$ to certain decay exponents that characterise the overall long-range connectivity.

\begin{theorem}[Finiteness vs.\ infiniteness of the percolation threshold]\label{thm:main_thm}
	Consider the weight-dependent random connection model $(\G_\beta\colon\beta>0)$ as given above.
	\begin{itemize}
		\item[(a)] The percolation threshold $\beta_c$ is finite, whenever 
		\[
		-\lim_{\mu\downarrow 0}\liminf_{n\to\infty}\frac{\log\left(\int_{n^{\mu-1}}^1\int_{n^{\mu-1}}^1 \rho(g(s,t)n) \,\textup{d}s\,\textup{d}t\right)}{\log n}<2.
		\]
		\item[(b)] The percolation threshold $\beta_c$ is infinite, whenever 
		\[
		-\lim_{\mu\uparrow 0}\limsup_{n\to\infty}\frac{\log\left(\int_{n^{\mu-1}}^1\int_{n^{\mu-1}}^1 \rho(g(s,t)n) \,\textup{d}s\,\textup{d}t\right)}{\log n}>2.
		\]
	\end{itemize}
\end{theorem}

We address the details of the proof of Theorem~\ref{thm:main_thm} in Section~\ref{secProofMainTheorem}. In fact, we prove the results under more general assumptions, see Propositions~\ref{thmInfinite} and~\ref{thmFiniteComponents}. First, we discuss some consequences and special instances of the result.

\begin{remark} \label{rem:MainResult} ~\
	\begin{enumerate}[(i)]
		\item 
		For many concrete choices of kernels and profiles discussed in the literature, the integrals appearing in the theorem are continuous in \(\mu\) and, moreover, both limits coincide. If this is the case, we define
		\begin{equation}
			\delta_\textup{eff} := -\lim_{n\to\infty}\frac{\log\left(\int_{1/n}^1\int_{1/n}^1 \rho(g(s,t)n) \,\textup{d}s\,\textup{d}t\right)}{\log n} \label{eqDeltaEff}
		\end{equation}
		and call \(\delta_\textup{eff}\) the \emph{effective decay exponent} associated with $\rho$ and $g$. Theorem~\ref{thm:main_thm} then simplifies to the two implications \(\beta_c<\infty\) if \(\delta_{\textup{eff}}<2\), and \(\beta_c=\infty\) if \(\delta_\textup{eff}>2\). In particular, continuous kernels and profiles of the form~\eqref{eq:rhodecay} lead to a well-defined \(\delta_\textup{eff}\). For profile functions of the form~\eqref{eq:rhodecay}, we further always have $\delta\geq \delta_{\textup{eff}}$, i.e., the inclusion of weights helps to create long-range connections and for strict inequality it is necessary that $g$ vanishes sufficiently fast at $(0,0)$.
		\item 
		The order (in \(n\)) of the integrals appearing in the theorem is independent of \(\beta\) for \(\rho\) satisfying~\eqref{eq:rhodecay}, except for the case when \(g(s,s)\asymp s\) as \(s\downarrow 0\) and \(\rho=\mathbbm{1}_{[0,1]}\) (or any other \(\rho\) of bounded support). In this situation we have for sufficiently large \(\beta\)
		\begin{equation*}
			\begin{aligned}
				\int\limits_{n^{-1}}^1\d s \int\limits_{n^{-1}}^1 \d t \ \1_{[0,1]}\big(\beta^{-1}g(s,t)n\big) \asymp \int\limits_{n^{-1}}^{\beta n^{-1}}\d s \int\limits_{n^{-1}}^{\beta n^{-1}} \d t \ 1 = (\beta-1)^2n^{-2}
			\end{aligned}
		\end{equation*}	 
		and hence \(\delta_\text{eff}=2\), thus this case is not covered by our results. 
		\item 
		The monotonicity conditions on $\rho$ and $g$ can be relaxed by straightforward domination arguments: for the existence of a supercritical phase, it suffices that $\rho\geq \rho_0, g\leq g_0$ with $\rho_0$ and $g_0$ satisfy the assumptions of Theorem~\ref{thm:main_thm}, and analogous but reverse bounds are sufficient for the absence of a supercritical phase.
		\item  
		{ Instead of working with a unit intensity Poisson process and varying \(\beta\), one could also vary the Poisson intensity \(\lambda>0\) and work with a fixed \(\beta\). The same proofs apply mutatis mutandis and we obtain \(\lambda_c<\infty\) if \(\delta_{\textup{eff}}<2\) and \(\lambda_c=\infty\) if \(\delta_{\textup{eff}}>2\). The same applies also to the refined results below in Propositions~\ref{thmInfinite} and~\ref{thmFiniteComponents}, in which the roles of \(\beta\) and \(\lambda\) can be exchanged in a similar manner.  }
	\end{enumerate}
\end{remark}

Let us give an intuitive explanation of the occurrence of \(\delta_{\textup{eff}}\). As \(n\) grows large, the minimum of \(n\) independent uniform random variables is roughly \(\nicefrac{1}{n}\) and consequently the double integral appearing in \eqref{eqDeltaEff} is essentially the probability of two randomly picked vertices from vertex sets of size $n$ at distance roughly \(n\) being connected by an edge. Ignoring additional correlations between edges arising from the introduction of vertex marks, the number of edges between the two vertex sets of size $n$ is roughly given by a binomial with \(n^2\) trials and connection probability \(n^{-\delta_\textup{eff}}\). If \(\delta_\textup{eff}<2\) the probability of having an edge between the two sets increases with $n$ whereas for \(\delta_\textup{eff}>2\) this probability decreases. The effective decay exponent \(\delta_{\textup{eff}}\) hence measures the occurrence of long edges in a way comparable to classical long-range percolation models, seen from a coarse-grained perspective. The truncation of the integral bounds in~\eqref{eqDeltaEff} is crucial to control the correlations arising from the vertex marks which is a necessity to identify the phase transition for the existence of a supercritical phase correctly. Indeed, at first glance, one might only want to calculate the decay exponent of the marginal distributions of single edges which is the rate at which the annealed probability of two typical vertices at distance \(n\) being connected decays. That exponent is
\begin{equation}
	\delta_{\textup{marg}}:=- \lim_{n\to\infty}\frac{\log \int_0^1 \d s \int_0^1 \d t \ \rho\big(g(s,t)n\big)}{\log n}. \label{eq:fakeexponent2}
\end{equation}
In analogy with homogeneous long-range percolation one might now assume that the graph contains an infinite cluster if \(\delta_{\text{marg}}<2\) and does not contain an infinite cluster if \(\delta_{\text{marg}}>2\). However, this does not take the aforementioned correlations into account and hence does not capture the behaviour accurately. By construction, we have \(\delta_{\text{marg}}\leq\delta_{\textup{eff}}\) and hence the correlations between edges make it harder to build an infinite component than in the independent case. Below we give an example in which we have strict inequality, showing that \(\delta_\text{marg}\) is indeed insufficient, see \eqref{eqDeltaIndBool}. It is an interesting observation that this is contrary to the physics literature mantra that positive correlations should help facilitate percolation, see e.g.\ \cite{PrakashEtAll1992} for some simulations of a correlated long-range percolation model in two dimensions, \cite{köhlerschindler2024criticalprobabilitiespositivelyassociated} for finite range models on trees, or \cite{DrewitzEtAll2022} for a proof for percolation induced by the Gaussian free field level sets on supercritical Galton-Watson trees. However, it is important to note that in the latter models the dependency structure is rather different from our model and also how the underlying geometry enters the construction differs considerably from our approach. 
\subsubsection*{The special case $\boldsymbol{\delta_{\textup{eff}}=2}$.}
The case $\delta_{\textup{eff}}=2$ is a generalisation of $\operatorname{LRP}(2)$ and the `$\nicefrac{1}{|x-y|^2}$-model' of long-range percolation \cite{AizenmanNewman1986}; we call this case \emph{pseudo-scale-invariant} based on the notion of scale-invariant long-range percolation. Indeed, if $g\equiv 1$ and $\rho(z)=1\wedge z^{-\delta}$, then $\delta_{\textup{eff}}=\delta$ and for $\delta=2$ we obtain a variant of $\operatorname{LRP}(2)$.

We do not attempt to cover the pseudo-scale-invariant case with the methods presented here, as it is clear that the existence of a percolation phase in this regime cannot be decided by looking at decay exponents alone: it is well known that the $\nicefrac{1}{|x-y|^2}$-model percolates for sufficiently large $\beta$ \cite{AizenmanNewman1986,DCGT2019} while on the other hand, a straightforward first moment argument, cf.\ \cite{Schulman1983}, yields that for any $g$ bounded away from $0$ the integrability condition
\[
\int_0^\infty\int_0^\infty \rho(|x-y|)\,\d x \,\d y <\infty
\]
alone is sufficient to imply $\beta_c=\infty.$ Hence choosing, say, $\rho(z)=1\wedge(z^{-2}\log(1+z)^{-2})$ yields $(\G_\beta\colon\beta>0)$ with $\beta_c=\infty$, whilst $\delta_{\textup{eff}}=\delta=2.$

In the examples discussed below, the $\delta_{\textup{eff}}=2$ regime usually corresponds to a boundary case, but there is at least one interesting example of a non-trivial kernel which yields models in the pseudo-scale-invariant regime for a large parameter range. This kernel corresponds to the age-dependent random connection model of \cite{GGLM2019}. It is an interesting open problem to identify the precise regime in which $\beta_c<\infty$ also for this kernel, but unfortunately there seems to be no easy way to adapt the method used in the present paper to cover model instances of interest in the pseudo-scale-invariant regime, even though the general approach \emph{does} work for the scale-invariant homogeneous model $\operatorname{LRP}(2)$, see \cite{DCGT2019}.

\subsubsection*{Analysis of concrete models from the literature.}
We introduce a new kernel that represents many of the established models for the right parameter choices and call it the \emph{interpolation kernel}. It is defined for \(\gamma\in[0,1)\) and \(\alpha\in[0,2-\gamma)\) and given by
\begin{equation}
	g_{\gamma,\alpha}(s,t) := (s\wedge t)^\gamma (s\vee t)^\alpha.  \label{eq:InterPol}
\end{equation}
In Table~\ref{tabInterPol}, we give various choices for \(\gamma\) and \(\alpha\) and the models they represent, depending on the geometric restrictions. With a similar calculation as in~\cite{GGLM2019} one can show that models constructed with the interpolation kernel have heavy-tailed degree distribution with power-law exponent \(\tau= 1 +(\nicefrac{1}{\gamma} \wedge \nicefrac{1}{(\alpha+\gamma-1)^+})\).

Let us calculate \(\delta_\textup{eff}\) for this kernel. Since \(g_{\gamma,\alpha}\leq 1\), we only consider profile functions as in \eqref{eq:rhodecay} with \(\delta>2\). For \(\delta\leq 2\), we have \(\beta_c<\infty\) by a simple coupling argument with homogenous continuum long-range percolation.  

\begin{table}[t!]
	\begin{center}
		\caption{Various choices for \(\gamma\), \(\alpha\) and \(\delta\) and the models they represent. Here, to shorten notation, \(\delta=\infty\) represents profile functions of bounded support. We do not distinguish whether the models were originally constructed on the lattice or a Poisson process.}
		\begin{tabular}{l r l}
			\toprule
			\textbf{Parameters} & \(\boldsymbol{\delta_\textup{eff}>2}\) \textbf{if} &\textbf{Names and references}    
			\tabularnewline
			\midrule
			\(\gamma = 0, \alpha = 0, \delta = \infty\) & always & random geometric graph \cite{Penrose1991,Penrose2003},
			\\
			& & Gilbert graph \cite{Gilbert61} 
			\vspace{4 pt}
			\tabularnewline
			\(\gamma = 0, \alpha =0, \delta<\infty\) & \(\delta>2\) & random connection model \cite{MeesterPenroseSarkar1996, Penrose2016}, 
			\\ 
			& & long-range percolation \cite{Schulman1983}
			\vspace{4pt}
			\tabularnewline
			\(\gamma>0, \alpha=0,\delta=\infty\)  & \(\gamma<1\) & Boolean model \cite{Hall85, Gouere08}, 
			\\
			& & scale-free Gilbert graph \cite{Hirsch2017}
			\vspace{4pt}
			\tabularnewline
			\(\gamma>0, \alpha = 0, \delta<\infty\) & \(\delta>2\), \(\gamma<\tfrac{\delta-1}{\delta}\) &soft Boolean model \cite{GGM21} 
			\vspace{4 pt}
			\tabularnewline
			\(\gamma = 0, \alpha>1, \delta = \infty\) & never &ultra-small scale-free geometric network \cite{Yukich2006} 
			\vspace{4 pt}
			\tabularnewline
			\(\gamma>0, \alpha=\gamma, \delta\leq \infty\) &  \(\delta>2\), \(\gamma<\nicefrac{1}{2}\) & scale-free percolation \cite{DeijfenHofstadHooghiemstra2013, DeprezWuthrich2019}, 
			\\ 
			& &geometric inhomogeneous random graphs \cite{BringmannKeuschLengler2019} 
			\vspace{4pt}
			\tabularnewline
			\(\gamma>0, \alpha=1-\gamma, \delta\leq\infty\) &  never & age-dependent random connection model \cite{GGLM2019}
			\tabularnewline
			\bottomrule
		\end{tabular}
		\label{tabInterPol}
	\end{center}
\end{table}


\begin{proposition} \label{lem:InterPol}
	For the interpolation kernel \(g_{\gamma,\alpha}\), with \(\gamma\in[0,1)\) and \(\alpha\in[0,2-\gamma)\), and a profile function \(\rho\) fulfilling~\eqref{eq:rhodecay} for some \(\delta>2\), we have 
	\[ 
	\delta_\textup{eff} = 
	\begin{cases}
		\delta, & \gamma<\tfrac{1}{\delta} \text{ and } \alpha<\tfrac{2}{\delta}-\gamma, 
		\\
		\delta(1-\gamma)+1, & \gamma > \tfrac{1}{\delta} \text{ and } \alpha<\tfrac{1}{\delta}, \\
		\delta(1-\gamma-\alpha)+2, & \alpha,\gamma>\tfrac{1}{\delta} \text{ or } \gamma<\tfrac{1}{\delta} \text{ and } \alpha>\tfrac{2}{\delta}-\gamma.
	\end{cases}
	\]
	In Particular, 
	\begin{enumerate}[(a)]
		\item If \(\gamma < 1-\nicefrac{1}{\delta}\) and
		\begin{enumerate}[(i)]
			\item if \(\alpha<1-\gamma\), then \(\delta_\textup{eff}>2\),
			\item if \(\alpha=1-\gamma\), then \(\delta_\textup{eff}=2\),
			\item if \(\alpha>1-\gamma\), then \(\delta_\textup{eff}<2\).
		\end{enumerate}
		\item If \(\gamma = 1- \nicefrac{1}{\delta}\) and \(\alpha\leq 1-\gamma\), then \(\delta_{\textup{eff}}=2\).
		\item If \(\gamma>1-\nicefrac{1}{\delta}\), then \(\delta_\textup{eff}<2\).
	\end{enumerate}	
\end{proposition}
\begin{proof}
	{  Pick a profile function of type~\eqref{eq:rhodecay}, note the symmetry of the kernel~\eqref{eq:InterPol}, and thus consider the integral
		\[
		\begin{cases}
			n^{-\delta} \int\limits_{n^{-1}}^1 \d s \int\limits_s^1 \d t \, s^{-\gamma\delta} t^{-\alpha\delta}, &  \text{ if }\gamma+\alpha \leq 1, \\
			\int\limits_{n^{-1}}^{n^{-{1}/{(\alpha+\gamma)}}} \d s \ n^{-1/\alpha}s^{-\gamma/\alpha} + n^{-\delta} \int\limits_{n^{-1}}^1 \d s \int\limits_s^1 \d t \, s^{-\gamma\delta} t^{-\alpha\delta}, &  \text{ if }\gamma+\alpha > 1.
		\end{cases} 
		\]
		Straightforward integration then yields the desired result by identifying the leading order terms in the different parameter regimes given in the proposition.   
	}    
\end{proof}

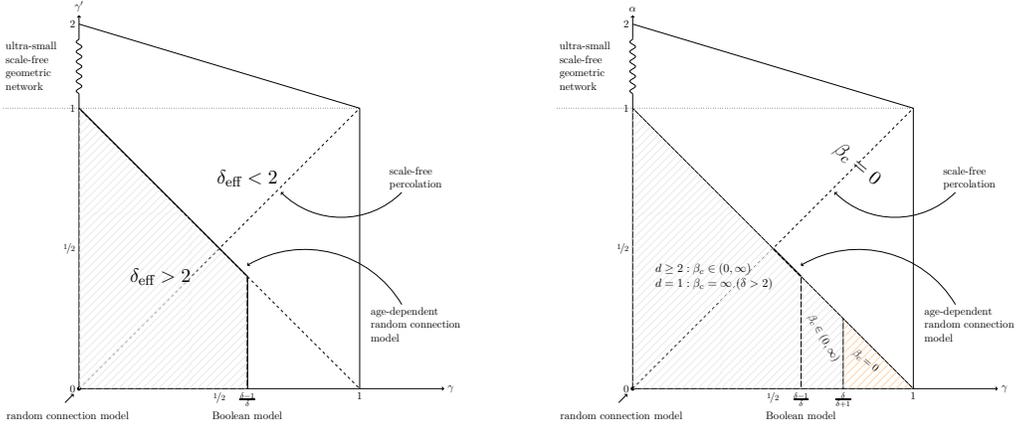
\begin{figure}[t!]
	\begin{center}
		\begin{subfigure}[t]{0.45\textwidth}
			\resizebox{\textwidth}{!}{
				\begin{tikzpicture}[every node/.style={scale=1}]
					\draw[->] (0,0) to (13,0) node[right] {$\gamma$};
					\draw	(5,0) node[anchor=north] {\nicefrac{1}{2}}
					(10,0) node[anchor=north] {1};
					
					\draw[dotted] (-2.7,10) to (10,10);
					\draw[] (10,0) to (10,10)
					(10,10) to (0,13);
					
					\draw[](0,0) to (0,10.5);
					\draw [->] (0,12.5) to (0,13.3) node[above] {$\alpha$};
					\draw[decorate, decoration = {snake, segment length = 10 pt, amplitude = 1mm}] (0, 10.5)--(0,12.5);
					\draw	(0,0) node[anchor=east] {0}
					(0,5) node[anchor=east] {\nicefrac{1}{2}}
					(0,10) node[anchor=east] {1}
					(0,13) node[anchor = east] {2};
					
					\draw[pattern=north east lines, pattern color=lightgray!30!, draw=none] 
					(0,0) plot[smooth,samples=200,domain=0:6,variable=\g]({0},{10*\g/6}) -- 
					plot[smooth,samples=200,domain=6:0,variable=\g]({6},{4*\g/6});
					
					\draw (0,8) to (4,4);
					\draw (4,4) to (4,0);
					\draw(4,4) to (10,4);
					
					\draw[dotted, thick, color = blue!100!black] (0,10) to (6,4);
					\draw[dashed] (6,4) to (10,0);
					\draw[dashed] (0,0) to (10,10);
					\draw (6,-0.7) node[align = left, anchor = north] {Boolean model};
					\draw (-1.7,11.5) node[align = left] {ultra-small \\ scale-free \\ geometric \\ network};
					\draw (12,2.3) node[align = left] {age-dependent \\ random connection \\ model};
					\draw[->, bend angle = 45, bend right] (11.5,3) to (6,4.4);
					\draw (12,7.5) node[align = left] {scale-free \\ percolation};
					\draw[->, bend angle = 45, bend left] (11.5,7) to (7.2,7);
					\draw (0,0) node[circle, fill = black, scale=0.4] {};
					\draw (-0.4,-0.7) node[align = left, anchor = north] {random connection model};
					\draw[->] (-0.5, -0.5) to (-0.2,-0.2);
					
					\draw (6, 0) node[anchor = north] {$1-\tfrac{1}{\delta}$};
					\draw[dotted, thick, color = blue!100!black] (6,0) to (6,4);
					
					\draw (4, 0) node[anchor = north] {$\tfrac{1}{\delta}$};
					\draw (0, 4) node[anchor = east] {$\tfrac{1}{\delta}$};
					\draw (0, 8) node[anchor = east] {$\tfrac{2}{\delta}$}; 
					
					\draw (2,3.5) node[scale = 1.5, thick]    {$\delta_\text{eff}=\delta$};
					\draw (4,8) node[scale = 1.5, thick]            {$\delta_\text{eff}=\delta(1-\gamma-\alpha)+2$};
					
					\draw (7,2) node[scale = 1.5, thick]    {$\delta_\text{eff}=\delta(1-\gamma)+1$};   
				\end{tikzpicture}
			}
			\caption{The three different values of \(\delta_{\textup{eff}}\) for the interpolation kernel subject to Proposition~\ref{lem:InterPol} for a profile function with \(\delta>2\). Shaded in gray the \(\delta_\textup{eff}>2\) phase of the model; the blue dotted line marks the pseudo-invariant regime \(\delta_\text{eff}=2\).}
			\label{subfig:deltaEff}
		\end{subfigure}
		\hfill
		\begin{subfigure}[t]{0.45\textwidth}
			\resizebox{\textwidth}{!}{  
				\begin{tikzpicture}[every node/.style={scale=1}]
					\draw[->] (0,0) to (13,0) node[right] {$\gamma$};
					\draw	(5,0) node[anchor=north] {\nicefrac{1}{2}}
					(10,0) node[anchor=north] {1};
					
					\draw[dotted] (-2.7,10) to (10,10);
					\draw[] (10,0) to (10,10)
					(10,10) to (0,13);
					
					\draw[](0,0) to (0,10.5);
					\draw [->] (0,12.5) to (0,13.3) node[above] {$\alpha$};
					\draw[decorate, decoration = {snake, segment length = 10 pt, amplitude = 1mm}] (0, 10.5)--(0,12.5);
					\draw	(0,0) node[anchor=east] {0}
					(0,5) node[anchor=east] {\nicefrac{1}{2}}
					(0,10) node[anchor=east] {1}
					(0,13) node[anchor = east] {2};
					
					\draw[thick] (0,10) to (10,0);
					\draw[dashed] (0,0) to (10,10);
					\draw (6,-0.7) node[align = left, anchor = north] {Boolean model};
					\draw (-1.7,11.5) node[align = left] {ultra-small \\ scale-free \\ geometric \\ network};
					\draw (12,2.3) node[align = left] {age-dependent \\ random connection \\ model};
					\draw[->, bend angle = 45, bend right] (11.5,3) to (6,4.4);
					\draw (12,7.5) node[align = left] {scale-free \\ percolation};
					\draw[->, bend angle = 45, bend left] (11.5,7) to (7.2,7);
					\draw (0,0) node[circle, fill = black, scale=0.4] {};
					\draw (-0.4,-0.7) node[align = left, anchor = north] {random connection model};
					\draw[->] (-0.5, -0.5) to (-0.2,-0.2);
					
					\draw (6, 0) node[anchor = north] {$\tfrac{\delta-1}{\delta}$};
					\draw[] (6,0) to (6,4);
					\draw (7.5, 0) node[anchor = north] {$\tfrac{\delta}{\delta+1}$};
					\draw[] (7.5,0) to (7.5,2.5);
					
					\draw[pattern=north east lines, pattern color=lightgray!30!, draw=none] 
					(0,0) plot[smooth,samples=200,domain=0:7.5,variable=\g]({0},{10*\g/7.5}) -- plot[smooth,samples=200,domain=7.5:0,variable=\g]({7.5},{2.5*\g/7.5});
					
					\draw[pattern=north west lines, pattern color=orange!30!, draw=none] 
					(0,0) plot[smooth,samples=200,domain=10:7.5,variable=\g]({10},{0}) -- plot[smooth,samples=200,domain=7.5:10,variable=\g]({7.5},{(\g-7.5)});
					
					\draw (8,8) node[scale = 1.8, thick, rotate = 320] {$\beta_c=0$};
					\draw (2.9,4) node[scale = 1.1, thick, align = left] {$d\geq 2: \beta_c\in(0,\infty)$ \\ $d=1: \beta_c=\infty \ (\delta>2)$};
					\draw (6.8,1.8) node[scale = 1, thick, rotate = 300] {$\beta_c\in(0,\infty)$};
					\draw (8.3,1) node[scale = 1, thick, rotate = 320] {$\beta_c=0$};
				\end{tikzpicture}
			}
			\caption{The regions where \(\beta_c\) is infinite, finite and positive, or zero, subject to Theorem~\ref{propPercolation}. The regions of Part~(a) are shaded, where orange (NW) marks Item~(i), grey (NE) marks Item~(ii). The \(d=1\) specifics only occur for \(\delta>2\).}
			\label{subfig:betaC}
		\end{subfigure}
	\end{center}
	\caption{Phase diagrams for the interpolation kernel depending on the values of \(\gamma\) and \(\alpha\). Dotted or dashed lines represent no change of behaviour. }
	\label{figPercolation}
\end{figure}

It is straightforward to extend the model to higher dimensions \(d\geq 2\): the vertex locations \(\eta\) are now given by a unit intensity Poisson process on \(\R^d\), and given the vertices and their i.i.d.\ marks, we connect any pair \((x,t_x)\) and \((y,t_y)\) of vertices independently with probability
\[
\rho\big(\beta^{-1}g_{\gamma,\alpha}(t_x,t_y)|x-y|^d\big).	
\]
In \cite{GLM2021} by Gracar et al., \cite{DeprezWuthrich2019} by Deprez and Wüthrich and \cite{Yukich2006} by Yukich, the question of the existence of subcritical percolation phases for certain special cases of our model is addressed. A perusal of their arguments shows, that they can easily be extended to the whole parameter regime of the interpolation kernel. Hence, we can give a rather comprehensive answer to the question when the interplay between geometric restrictions and additional edges coming from the marks lead to a non-trivial phase transition, see Figure~\ref{subfig:betaC}.

\begin{theorem} \label{propPercolation}
	Let \((\G_\beta\colon\beta>0)\) be the weight-dependent random connection model constructed with a profile function \(\rho\) fulfilling~\eqref{eq:rhodecay} for some \(\delta>1\) and the interpolation kernel \(g_{\gamma,\alpha}\) for \(\gamma\in[0,1)\) and \(\alpha\in[0,2-\gamma)\). 
	\begin{enumerate}[(a)]
		\item For \(\alpha\leq 1-\gamma\), we have:
		\begin{enumerate}[(i)]
			\item If \(\gamma>\tfrac{\delta}{\delta+1}\), then \(\beta_c=0\). 
			\item If \(\gamma<\tfrac{\delta}{\delta+1}\) and
			\begin{itemize}
				\item either \(d\geq 2\), or \(d=1\) and \(\delta<2\), or \(d=1\) and \(\gamma>1-\nicefrac{1}{\delta}\), or \(d=1\) and \(\alpha=1-\gamma\) and \(\gamma\geq \nicefrac{1}{2}\), then \(\beta_c\in(0,\infty)\);
				\item \(d=1\) and \(\alpha<1-\gamma\) and \(\delta>2\) and \(\gamma<1-\nicefrac{1}{\delta}\), then \(\beta_c=\infty\).
			\end{itemize}
		\end{enumerate}
		\item For \(\alpha>1-\gamma\), we have \(\beta_c=0\). 
	\end{enumerate}
\end{theorem}

\begin{proof}
	The proof of part (a)(i) is done in~\cite{GLM2021}, (a)(ii) is a combination of \cite{GLM2021}, known results about the existence of supercritical percolation phases in \(d\geq 2\)~\cite{MeesterRoy1996} and in long-range percolation in \(d=1\) when \(\delta\in(1,2)\)~\cite{NewmanSchulman1986}, our main Theorem~\ref{thm:main_thm}, and Corollary~\ref{cor:AgeDependent} below. Part (b) can be easily proven by generalising the proof for \(\alpha=\gamma>\nicefrac{1}{2}\) of~\cite{DeprezWuthrich2019} to all \(\alpha>1-\gamma\).
\end{proof}

By combining Proposition~\ref{lem:InterPol} and Theorem~\ref{propPercolation}, we can see that the robust regime \(\beta_c=0\) for the interpolation kernel model is always contained in the regime where \(\delta_\text{eff}<2\). Moreover, while the present manuscript underwent peer review, it was shown in \cite{jacob2025subcriticalannulus} that robustness implies $\delta_\textup{eff}\leq 2$ in any dimension for a large class of translation invariant inhomogeneous percolation models that includes the weight-dependent random connection model. We believe that there are in fact no robust locally finite instances of the weight-dependent random connection model in the pseudo-scale-invariant regime \(\delta_\text{eff}=2\), but currently have no proof of this conjecture.

To conclude this section, we further discuss two models of particular interest that can be represented by the interpolation kernel. Our focus is, again, on dimension \(d=1\) and profile functions satisfying~\eqref{eq:rhodecay} for \(\delta>2\).

\begin{figure}
	\begin{center}
		\resizebox{0.9\textwidth}{!}{
			\begin{tikzpicture}[scale=0.65, every node/.style={scale=0.5}, show background rectangle]
				\node (X) at (0,1.5)[circle, fill=blue, label = {$\x$}] {};
				\node (Y) at (3,3) [circle, fill=blue, label = {below:$\y$}] {};
				\node (Z) at(6.3, 1) [circle, fill=blue, label = {right: $\z$}] {};
				\draw[color = white] (X) circle (4);
				\draw[color = white] (Z) circle (4.2);
				\draw[] (X) circle (3.7);
				\draw[] (Y) circle (1.6);
				\draw[] (Z) circle (2);
				\draw[thick] (X) to (Y);
			\end{tikzpicture}
			\hspace{3 cm}
			\begin{tikzpicture}[scale=0.65, every node/.style={scale=0.5},show background rectangle]
				\node (X) at (0,1.5)[circle, fill=blue, label = {$\x$}] {};
				\node (Y) at (3,3) [circle, fill=blue, label = {below:$\y$}] {};
				\node (Z) at(6.3, 1) [circle, fill=blue, label = {right: $\z$}] {};
				\draw[color = white] (Z) circle (4.2);
				\draw[] (X) circle (3.7);
				\draw[] (Y) circle (1.6);
				\draw[] (Z) circle (2);
				\draw[thick] (X) to (Y);
				\draw[dashed] (X) circle (4);
				\draw[dashed] (X) to (3.9,0.577);
				\node (E1) at (2,0.5) {$R(x,z)t_x^{-\gamma/2}$};
			\end{tikzpicture}
			\hspace{3 cm}
			\begin{tikzpicture}[scale=0.65, every node/.style={scale=0.5},show background rectangle]
				\node (X) at (0,1.5)[circle, fill=blue, label = {$\x$}] {};
				\node (Y) at (3,3) [circle, fill=blue, label = {below:$\y$}] {};
				\node (Z) at(6.3, 1) [circle, fill=blue, label = {right: $\z$}] {};
				\draw[color = white] (X) circle (4);
				\draw[] (X) circle (3.7);
				\draw[] (Y) circle (1.6);
				\draw[] (Z) circle (2);
				\draw[thick] (X) to (Y);
				\draw[dashed] (Z) circle (4.2);
				\draw[dashed] (Z) to (2.2, 0);
				\node (E2) at (4.2,0.1) {$R(y,z)t_z^{-\gamma/2}$};
				\draw[thick] (Z) to (Y);
			\end{tikzpicture}
		}
	\end{center}
	\caption{Example for connection mechanism of the soft Boolean model in two dimensions. The solid lines represent the edges of the graph.}
	\label{fig:SoftBool}
\end{figure}
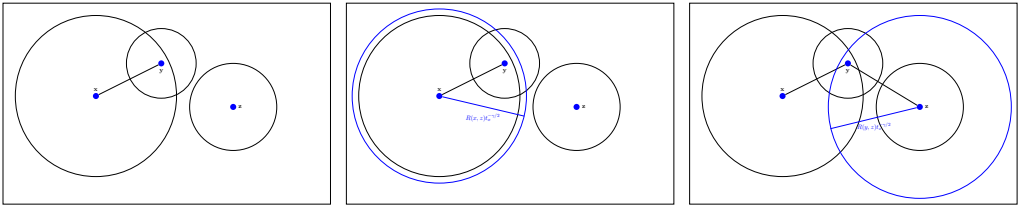

\noindent\emph{Soft Boolean model:} This model is a long-range variant of the classical (`hard') Poisson-Boolean model. It corresponds to the kernel \(g_{\gamma,0}(s,t)=(s\wedge t)^\gamma\) for some \(\gamma\in(0,1)\) and can be interpreted as follows. Each vertex \(x\) is assigned an independent radius \( t_x^{-\gamma/d}\). Given all vertices, each potential edge \(\{x,y\}\) is assigned an independent heavy-tailed random variable \(R(x,y)\) with \(\P\{R(x,y)>r\}=r^{-d\delta}\) (for \(r>1\)) and any given pair of vertices \(x,y\) is connected if the vertex with the smaller assigned radius, say \(x\), is contained in the ball of radius \(\beta^{1/d} R(x,y)t_y^{-\gamma/d}\) centred in \(y\). That is
\[
\x\sim\y \Leftrightarrow |x-y| \leq \beta^{1/d} R(x,y)\big(t_x^{-\gamma/d}\vee t_y^{-\gamma/d}\big),
\]

see Figure~\ref{fig:SoftBool}. A variant closer to the classical Poisson-Boolean model is given by the \emph{sum kernel} \(g^{\text{sum}}(s,t)=(s^{-\gamma/d}+t^{-\gamma/d})^{-d}\) where an edge is drawn if the ball with radius \(\beta^{1/d} R(x,y)t_x^{-\gamma/d}\) centred in \(x\) and the ball with radius \(\beta^{1/d} R(x,y) t_y^{-\gamma/d}\) centred in \(y\) intersect, cf.\ \cite{GGM21} for a more detailed discussion. As \(g^{\text{sum}}\leq g_{\gamma,0}\leq 2^d g^\text{sum}\) both kernels show qualitatively the same behaviour in the question of existence of a non trivial percolation phase transition. It follows from Theorem~\ref{propPercolation} that, for \(\delta>2\), the one-dimensional soft Boolean model provides three open parameter regimes: one where the graph always contains an infinite component, an intermediate regime where an infinite component exists for large values of $\beta$, and one where no infinite component can exist. Allowing \(\gamma=0\), this model also includes the homogeneous random connection model. 

Let us further calculate the decay exponent of the marginal distribution of single edges as mentioned in the explanation around~\eqref{eq:fakeexponent2}. If we choose \(\rho(x)\asymp 1\wedge x^{-\delta}\) for \(\delta>2\) and some \(\nicefrac{1}{2}<\gamma<1-\nicefrac{1}{\delta}\), we calculate 
\begin{equation}\label{eqDeltaIndBool}
	\begin{aligned}
		\int_0^1 \d s \int_0^1 \d t \ \big(1\wedge (s\wedge t)^{-\gamma\delta}n^{-\delta}\big) & \asymp n^{-\tfrac{1}{\gamma}} + n^{-\delta}\int_{n^{-1/\gamma}}^1 \d s \ s^{-\gamma\delta} \asymp n^{-\tfrac{1}{\gamma}}. 
	\end{aligned}
\end{equation}
Since \(\gamma>\nicefrac{1}{2}\), we have \(\delta_{\text{marg}}<2\) and we have \(\beta_c<\infty\) for the independent edges model where each radius is resampled independently for any potential edge, but as \(\gamma<1-\nicefrac{1}{\delta}\), we have \(\delta_{\textup{eff}}>2\) by Proposition~\ref{lem:InterPol} and hence \(\beta_c=\infty\) by Theorem~\ref{thm:main_thm} for the actual model where each radius is sampled just once.

\noindent\emph{Age-dependent random connection model:} This model was introduced in~\cite{GGLM2019} as the weak local limit of an age-based spatial preferential attachment model and corresponds to the choice \(\gamma\in(0,1)\) and \(\alpha=1-\gamma\). Here, the vertex marks play the role of vertex birth times and small marks correspond to early birth times and hence to old and thus influential vertices. By Theorem~\ref{propPercolation} we have \(\beta_c=0\) if \(\gamma>\nicefrac{\delta}{(\delta+1)}\) and \(\beta_c\in(0,\infty)\) if \smash{\(\gamma\in\big(\tfrac{\delta-1}{\delta},\tfrac{\delta}{\delta+1}\big)\)}. However, we have \(\delta_\textup{eff}=2\) for the whole parameter regime \(0<\gamma\leq 1-\nicefrac{1}{\delta}\). This illustrates a scale-invariance property that is built into the model and is a result of the dynamics within the model coming from the role of vertex marks being birth times. Since our main theorem does not capture the pseudo-scale-invariant regime, we can only deduce results by comparing the age-dependent random connection model to models where the finiteness or infiniteness of \(\beta_c\) is known. To this end, we consider the scale-free percolation model, introduced by Deijfen et al.~\cite{DeijfenHofstadHooghiemstra2013} also known later as infinite geometric inhomogeneous random graph~\cite{BringmannKeuschLengler2019},  which coincides with the choice of \(\gamma=\alpha\) in the interpolation kernel. Here, it is known that \(\beta_c=0\) if \(\gamma>\nicefrac{1}{2}\) and \(\beta_c=\infty\) if \(\gamma<\nicefrac{1}{2}\). Now note that scale-free percolation and the age-dependent random connection model coincide for \(\gamma=\nicefrac{1}{2}\). In \cite{KomjathyLodewijks2020, BringmannKeuschLengler2019,HofstadHoornMaitra2021} it is shown that the `KPKVB-model', a hyperbolic random graph model, has one-dimensional scale-free percolation as weak local limit (after a suitable change from hyperbolic to weighted Euclidean coordinates). Results from Bode et al.~\cite{BodeFountoulakisMuller2015} for the KPKVB-model then show that there exists an infinite cluster in the infinite geometric inhomogeneous random graph for \(\rho(r)=\mathbbm{1}_{[0,1]}(r)\), whenever \(\gamma=\nicefrac{1}{2}\) and \(\beta\) is sufficiently large. By monotonicity in \(\gamma\) of \(g_{\gamma,1-\gamma}\), it follows that the same holds true for the age-dependent random connection model for any \(\gamma\geq \nicefrac{1}{2}\) and any $\rho$ that is bounded away from $0$ in a neighbourhood of $0$. It remains an interesting open problem to show whether there can be percolation for \(\gamma<\nicefrac{1}{2}\) in this setting. 

The discussion of this paragraph is summarised in the following corollary.

\begin{corollary}\label{cor:AgeDependent}
	Let \((\G_\beta\colon\beta\geq 0)\) be the weight dependent random connection model in dimension \(d=1\), constructed with the interpolation kernel \(g_{\gamma,\alpha}\) and any non-trivial profile function \(\rho\) satisfying the upper bound in~\eqref{eq:rhodecay} for some \(\delta\in(2,\infty)\). If either \(\gamma=\alpha=\nicefrac{1}{2}\), or \(\nicefrac{1}{2}<\gamma<\nicefrac{\delta}{(\delta+1)}\) and \(\alpha=1-\gamma\), we have \(\beta_c\in(0,\infty)\). 
\end{corollary}

\subsubsection*{Giant components in finite versions of the model.}
Another application of our results is the existence of a component of linear size in finite versions of our model, also known as the \emph{giant component}. If a giant exists, one is interested in whether it is unique and how many vertices belong to it. The paper \cite{HofstadHoornMaitra2021} elaborates how weight-dependent random connection-type models arise as weak local limits of models on finite domains. There, growing sequences of graphs \((\G_\beta^{(n)}\colon n\in\N)\) are constructed, where \(\G_\beta^{(n)}\) consists of \(n\) vertices which are independently placed into the unit interval \((-\nicefrac{1}{2},\nicefrac{1}{2})\). Translated to our parametrisation each vertex carries an independent uniform mark and, given locations and marks, each pair \((x,t),(y,s)\) of vertices is connected independently with probability \(\rho(\beta^{-1}g(s,t)n|x-y|)\). The scaling factor \(n\) ensures that the graph remains sparse and the \(n\) vertices can hence be considered as being embedded into \((-\nicefrac{n}{2},\nicefrac{n}{2})\). Recall that \(\theta(\beta)\) denotes the percolation probability in $\G_\beta$, it is formally defined below in~\eqref{eq:percoProb}. Our proofs yield the following corollary regarding finite versions of the weight-dependent random connection model.

\begin{corollary}\label{cor:finiteGraphs}
	Let \((\G_\beta^{(n)}\colon n\in\N)\) be the above sequence of finite graphs on intervals with weak local limit $\G_\beta$ given by an instance of the weight-dependent random connection model with kernel $g$ and profile $\rho.$ If $\rho$ and $g$ satisfy the assumptions of Theorem~\ref{thm:main_thm}(a), then \((\G_\beta^{(n)}\colon n\in\N)\) contains a connected component that grows linearly in \(n\) for large enough \(\beta\). Conversely, if $\rho$ and $g$ satisfy the assumptions of Theorem~\ref{thm:main_thm}(b), then there is no giant component in the above sense for any $\beta>0.$
\end{corollary}

It is well-known (see for example \cite{vdH2021Giant}) that a graph sequence cannot have a giant component if its weak local limit does not percolate, hence it remains to prove the statement in the \(\delta_\textup{eff}<2\) case, which is done at the end of Section~\ref{secInfiniteComponent}.

\subsubsection*{Organisation of the remainder of the paper.} 
In Section~\ref{SecGen}, we start with a precise mathematical construction of the model and introduce more general vertex sets. We then state in Section~\ref{secProofMainTheorem} two propositions: a sufficient condition for the existence of a supercritical phase, Proposition~\ref{thmInfinite}, and a sufficient condition for the non-existence of such a phase, Proposition~\ref{thmFiniteComponents}. Combined, these imply Theorem~\ref{thm:main_thm}. We prove Proposition~\ref{thmInfinite} in Section~\ref{secInfiniteComponent} and Proposition~\ref{thmFiniteComponents} in Section~\ref{secFiniteComponent}. Let us remark that we give the proof of Proposition~\ref{thmInfinite} before the proof of Proposition~\ref{thmFiniteComponents} as it contains the main new contributions of the paper; both mathematically and from the point of view of specific model instances, particularly through Lemma~\ref{lemLongConnection} and Lemma~\ref{lem:Decorrelate}. However, the proof of absence of percolation is conceptionally easier and the reader may want to start with Section~\ref{secFiniteComponent} if they want to familiarise themselves with the model first.


\section{General set up and proof of {Theorem~\ref{thm:main_thm}}} \label{SecGen}
We now introduce our model in a rigorous way and define the properties we require of the underlying point process for our results to hold. 
We then formulate Propositions~\ref{thmInfinite} and~\ref{thmFiniteComponents}, which imply Theorem~\ref{thm:main_thm}. 

\subsection{Formal construction of the model} \label{sec:formalConstr}
Let $\eta$ be a standard renewal process on $\R$ with intensity $\lambda\in(0,\infty)$. In the following, as the intensity plays no particular role beyond being positive and finite, we omit it from our notation.
We write \(\eta_0\) for the \emph{Palm version} of $\eta$ containing a point at $0\in\R$, see e.g.\ \cite{Daley2003} for background on Palm distributions. { That is, $\eta_0=\{X_j\colon j\in\Z$ such that $X_k<X_\ell$ for $k<\ell\}$ with $X_0=0$ and the differences \(X_{j}-X_{j-1}\) are independent and identically distributed. In the literature these differences are often referred to as innovations or holding times; we call them differences as this aligns best with the interpretation of distances between vertices. We think of the vertex in location $0$ as the `root vertex' of its cluster.}
The weight-dependent random connection model is constructed as a deterministic functional $\G_{\beta,\rho,g}$ of the points of $\eta_0$, and two independent i.i.d.\ sequences of edge and vertex marks as follows. Let $\mathscr{T}=\{T_j\colon j\in \Z\}$ be a family of i.i.d.\ random variables distributed uniformly on $(0,1)$ independent of $\eta_0$ and define the marked point process
\[
\mathscr{X}_0:=\{\X_j=(X_j,T_j)\in\eta_0\times\mathscr{T}, j\in\Z \}.
\]
Hence, $\mathscr{X}_0$ is a point process on $\R\times(0,1)$ with intensity measure $\lambda\otimes\operatorname{Uniform(0,1)}$ and a distinguished, `typical' point located at the origin. Let further $\mathscr{U}=\{U_{i,j}\colon i<j\in\Z\}$ be a second family of i.i.d.\ Uniform$(0,1)$ random variables, independent of $\mathscr{X}_0$, that we call \emph{edge marks}. We use \(\mathscr{X}_0\) and \(\mathscr{U}\) to define a point process 
\begin{equation}
	\xi_0:= \Big\{\big(\{\X_i,\X_j\},U_{i,j}\big)\in\mathscr{X}_0^{[2]}\times\mathscr{U}\colon i<j\in\Z\Big\}, \label{eq:edgeMarking}
\end{equation}
where \(\mathscr{X}_0^{[2]}\) denotes the set of all subsets of size two of \(\mathscr{X}_0\). We call \(\xi_0\) an \emph{independent vertex-edge-marking} of \(\eta_0\) in accordance with~\cite{HvdHLM20}. Observe that \(\mathscr{X}_0\) as well as \(\eta_0\) can be recovered from $\xi_0$. From here on onwards, we assume that we work on a probability space on which the vertex-edge-marking \(\xi_0\) is defined and denote the underlying probability measure by \(\P\) and the corresponding expectation by \(\E\).

Now fix $\beta>0$, a profile function $\rho$ and a kernel function $g$. Then $\G_{\beta,\rho,g}(\xi_0)$ is the graph with vertex set $\mathscr{X}_0$ and edge set
\[\Big\{\{\X_i,\X_j\}\colon U_{i,j}\leq \rho\big(\tfrac{1}{\beta}\,g(T_i,T_j)\,|X_i-X_j|\big), i<j\Big\}.\] 
To keep the notation concise, we write $\cG_\beta=\cG_{\beta,\rho,g}(\xi_0)$. Note that this graph has the same law as the previous one if \(\eta\) is a Poisson process, justifying the slight abuse of notation. 
Due to our a priori ordering of the vertices by their locations, one may think of \(\mathscr{X}_0\) as a marked lattice by considering the vertex indices instead of their actual locations on the real line. Our proofs therefore apply immediately to $\Z$-based versions of the model as well. {Indeed, the lattice can formally be viewed as a renewal process with all inter-point distances equal to one (or i.i.d.\ geometrically distributed in case of a Bernoulli-site percolated lattice).}  However, note that even though we define the model using a specific ordering of $\eta_0$, the distribution of $\cG_\beta$ is independent of the chosen order.

Our central requirement on the renewal process $\eta$ is the following regularity condition:

\begin{definition}[Evenly spaced renewal process]\label{DefEvenlySpaced} 
	Let $\eta$ be a standard renewal process on \(\R\) and denote by $\mathbf{P}_0$ the law of its Palm version \(\eta_0\). We say that $\eta$ (and also $\eta_0$) is \emph{evenly spaced}, if for any \(a>2\lambda\) the large-deviation bound
	\begin{equation} \label{eq:ES}
		\mathbf{P}_0\{|X_{-n}-X_{n-1}|>a n \} = o(n^{-2}) ,
	\end{equation}
	as \(n\to\infty\), holds. Similarly, we call a realisation of \(\eta_0=\omega=\{x_i\colon i\in\Z\}\), for which 
	\begin{equation} \label{eq:propspaced}
		|x_{-n}-x_{n-1}|< a n,
	\end{equation} 
	is satisfied for some \(a\), \(n\)\emph{-properly spaced}.
\end{definition}
\begin{remark} \label{rem:evenly}
	\begin{enumerate}[(i)]
		\item { Property~\eqref{eq:ES} essentially says that distances between vertices (i.e., the differences of consecutive locations) are not too large and do not fluctuate too much on large scales. This then guarantees that vertex locations are sufficiently dense uniformly over all scales, which is required to prove the existence of a supercritical phase.
			\item All renewal processes with light-tailed differences, including the Poisson process, are evenly spaced. In case of heavy-tailed differences of the subexponential class, the existence of their third moment is equivalent to the evenly spaced property~\cite{DenisovEtal_2008_ldp}.} 
		\item Theorem~\ref{thm:main_thm} holds true for the weight-dependent random connection model based on any evenly spaced renewal process, any deterministic point set satisfying the corresponding deterministic spacing conditions, { and any stationary simple point process that is appropriately stochastically dominated by an evenly spaced renewal process. An example of the latter could, for instance, be a Cox process, where the random intensity measure first generates a random closed set and then places the Poisson points with two different intensities for regions within said set and its complement; see e.g.~\cite{JahnelTobiasCali2022}.}
	\end{enumerate}	 
\end{remark}

\subsection{Strengthening Theorem~\ref{thm:main_thm}} \label{secProofMainTheorem}
We begin by formalising the notion of percolation. For two vertices $\X_i,\X_j\in\mathscr{X}_0$, we denote by $\{\X_i\sim \X_j\text{ in }\cG_\beta\}$ the event that $\X_i$ and $\X_j$ are connected by an edge in $\cG_\beta$. When the graph \(\cG_\beta\) is fixed, we simply write \(\{\X_i\sim\X_j\}\). We define $\{\0\leftrightarrow\infty \text{ in }\cG_\beta\}$ as the event that the root $\X_0$ is the starting point of an infinite self-avoiding path $(\X_0, \X_{i_1}, \X_{i_2},\dots)$ in $\cG_\beta$, i.e., $\X_{i_j}\neq\X_{i_k}$ for all $j\neq k$ and $\X_{i_{j}}\sim \X_{i_{j+1}}$ for all $j\geq 0$. We set
\begin{equation} \label{eq:percoProb}
	\theta(\beta) = \P\{\0\leftrightarrow\infty \text{ in }\cG_\beta\}, \quad \beta>0.
\end{equation}
If $\theta(\beta)>0$, then $\cG_\beta$ percolates almost surely by ergodicity. Conversely, if $\cG_\beta$ contains an infinite component, $\X_0$ is connected to it with positive probability and, hence, $\theta(\beta)>0$. We therefore call $\theta(\beta)$ the percolation probability and define the \emph{percolation threshold} \(\beta_c\) as 
\begin{equation*}
	\beta_c:=\inf\{\beta>0\colon\theta(\beta)>0\}. 
\end{equation*}	

We now state two propositions which combined prove a strengthened version of Theorem~\ref{thm:main_thm}. We start with a sufficient condition for the existence of an infinite cluster.

\begin{proposition}[Existence of supercritical phase]\label{thmInfinite}
	Let \(\eta\) be an evenly spaced renewal process in the sense of Definition~\ref{DefEvenlySpaced}. Let \(\xi_0\) be an independent vertex-edge-marking of $\eta$, \(\rho\) be a profile function and \(g\) be a kernel function. Let \(K_n:=(n!)^3 K^n, n\in\N\) for some \(K\in\N\). Assume that there exists \(\mu\in(0,\nicefrac{1}{2})\) such that
	\begin{equation}
		\sup_{n\geq 2} n^3 K\exp\Big(-K_{n-1}^2 \int\limits_{2K_{n-1}^{\mu-1}}^1 \d s \int\limits_{2K_{n-1}^{\mu-1}}^1 \d t \ \rho\big(g(s,t )K_n\big)\Big) =o(1), \text{ as }K\to\infty. \tag{A1} \label{eqInfinite}
	\end{equation}		
	Then, for \((\cG_{\beta, \rho, g}(\xi_0)\colon\beta>0)\), we have $\beta_c<\infty$.
\end{proposition}

The proof of Proposition~\ref{thmInfinite} is given in Section~\ref{secInfiniteComponent} and follows a strategy similar to that of the proof of \cite[Theorem~1(i)]{DCGT2019}. Condition \eqref{eqInfinite} can be seen as a generalised version of \cite[Equation~(8)]{DCGT2019} and precisely quantifies the order of the probability that certain long edges are absent in the graph that we need for our proof to work. The analogous result for the non-existence of an infinite component is given next.

\begin{proposition}[Absence of supercritical phase] \label{thmFiniteComponents}
	Let \(\eta\) be a standard renewal process on \(\R\) with corresponding vertex-edge-marking \(\xi_0\) and let \(g\) be a kernel function.
	Assume further, that \(\rho\) is a profile function with \(\rho(0+)<1\) and that there exists \(\mu\in(0,\nicefrac{1}{2})\) such that
	\begin{equation}
		\sum_{n\in\N} 		2^{2n}\int_{2^{-(1+\mu)n}}^1 \d s \int_{2^{-(1+\mu)n}}^1 \d t \ \rho(g(s,t)2^n) <\infty.	 \tag{A2} \label{eqFinite}
	\end{equation}
	Then for \((\cG_{\beta, \rho, g}(\xi_0)\colon\beta>0)\), we have \(\beta_c=\infty\).
\end{proposition}
\begin{remark}\label{rem:finiteComp}
	The proof of Proposition~\ref{thmFiniteComponents}, in fact, works for graphs constructed on any stationary and ergodic simple point process. While we stick here to the renewal process to keep notation light, we comment on that matter below the proof in Section~\ref{secFiniteComponent} in Remark~\ref{rem:ProofFiniteComp}.
\end{remark}

The assumption \(\rho(0+)<1\) in Proposition~\ref{thmFiniteComponents} is a technical requirement needed in the proof of Proposition~\ref{thmFiniteComponents} below and can essentially be viewed as a continuum version of the analogous requirement in long-range percolation on $\Z$ that not all nearest-neighbour-edges be present. If \(\eta\) is the Poisson point process the additional condition $\rho(0+)<1$ can be dropped.

\begin{corollary}\label{CorolFiniteCompPois} 
	Let \(\eta\) be a standard Poisson point process. Let \(\xi_0\) be an independent vertex-edge-marking based on $\eta$, \(g\) be a kernel and \(\rho\) be a profile function such that assumption~\eqref{eqFinite} of Proposition~\ref{thmFiniteComponents} is fulfilled. Then for \((\cG_{\beta,\rho,g}(\xi_0)\colon\beta>0)\), we have \(\beta_c=\infty\).
\end{corollary}

We are now able to prove the main Theorem~\ref{thm:main_thm} using the previous propositions. 

\begin{proof}[Proof of Theorem~\ref{thm:main_thm}]
	Let us define
	\[
	\delta_\text{eff}^+(\mu):=-\liminf_{n\to\infty}\frac{\log\left(\int_{n^{\mu-1}}^1\int_{n^{\mu-1}}^1 \rho(g(s,t)n) \,\textup{d}s\,\textup{d}t\right)}{\log n}
	\]
	and assume \(\delta_\text{eff}^+(0+)<2\). Observe that \(2K_{n-1}^{\mu-1}=2K_n^{\mu-1}(n^3K)^{1-\mu}\leq K_n^{\mu+\mu(K)-1}\) for some small  \(\mu(K)\downarrow 0\) as \(K\to\infty\). Choosing \(K\) large and \(\mu\) small enough such that \(\delta_\text{eff}^+(\mu+\mu(K))<2\), we have
	\[
	K_{n-1}^2 \int\limits_{2K_{n-1}^{\mu-1}}^1 \d s \int\limits_{2K_{n-1}^{\mu-1}}^1 \d t \ \rho\big(g(s,t )K_n\big) \geq c\,  K_{n-1}^{2-\delta_{\textup{eff}}^+(\mu+\mu(K))},
	\]
	uniformly in \(n\) for some constant \(c>0\). Hence Assumption~\eqref{eqInfinite} is satisfied and Theorem~\ref{thm:main_thm}(a) is a consequence of Proposition~\ref{thmInfinite}. 
	Similarly, define
	\[
	\delta_\textup{eff}^-(\mu):=-\limsup_{n\to\infty}\frac{\log\left(\int_{n^{-\mu-1}}^1\int_{n^{-\mu-1}}^1 \rho(g(s,t)n) \,\textup{d}s\,\textup{d}t\right)}{\log n}
	\]
	and assume \(\delta_\textup{eff}^-(0+)>2\). Then, for small enough \(\mu\) and large enough \(n\), we have for some constant \(C>0\)
	\[
	2^{2n}\int_{2^{-(1+\mu)n}}^1 \d s \int_{2^{-(1+\mu)n}}^1 \d t \ \rho(g(s,t)2^n)\leq C 2^{n(2-\delta_\text{eff}^-(\mu))}.
	\]
	Hence, Assumption~\ref{eqFinite} is satisfied and the proof of Theorem~\ref{thm:main_thm} is concluded by applying Corollary~\ref{CorolFiniteCompPois}.
\end{proof}

\section{Proof of Proposition \ref{thmInfinite}} \label{secInfiniteComponent}

We now proceed to prove the existence of an infinite component under the assumptions of Proposition~\ref{thmInfinite}. Our strategy is based on the approach developed in \cite[Theorem~1(i)]{DCGT2019} to establish the existence of supercritical percolation regime in classical long-range percolation LRP($\delta$) in the scale-invariant regime \(\delta=2\). The proof involves a multi-scale argument working roughly as follows: at stage $n$, the lattice is covered by half-way overlapping blocks of \(K_n\) lattice points. The overlap has the effect that if two adjacent blocks contain a rather dense connected component each, the two components necessarily intersect due to the pigeon hole principle. Since enough blocks at stage \(n\) have such components, a larger block of size \(K_{n+1}\) containing several \(K_n\)-blocks must also contain a component of positive (but ever so slightly smaller) density. If the loss of density at each stage can be kept sufficiently small, the desired result follows by iterating this construction and taking the limit \(n\to\infty\). However, renormalisation requires initialisation with local clusters of a given (large) density, hence one cannot hope to obtain quantitative bounds for $\beta_c$ using this technique.

In our model, we face the challenge of the additional randomness of the marks. We demonstrate below that a modified version of the above strategy works under the assumptions of Proposition~\ref{thmInfinite}, i.e., in particular if \(\delta_{\textup{eff}}<2\). Although we are able to control the influence of the marks by a carefully tailored argument, the error probabilities arising at each stage become too large if one moves into the pseudo-scale-invariant regime, hence the approach breaks down if $\delta_{\text{eff}}=2$.

For $N\in\N$ and $i\in \Z$ let
\[
\mathscr{B}_N^i:=\{\X_{N(i-1)},\dots,\X_{Ni},\dots,\X_{N(i+1)-1}\}
\]
and $\mathscr{B}_N:= \mathscr{B}_N^0 = \{\X_{-N},\dots, \X_{N-1}\}$. Each set \(\mathscr{B}_N^i\) consists of precisely \(2N\) consecutive vertices. If \(\eta\) is the lattice, then \(X_j=j\) for each \(j\in\Z\) and \(\mathscr{B}_N^i\) is simply the lattice interval \([N(i-1), N(i+1))\cap \Z\), matching the notation of \cite{DCGT2019}. In the general setting, the sets \(\mathscr{B}_N^i\) are blocks of vertices that all contain the same number of vertices but with random distances between consecutive vertices. Note that two consecutive blocks $\mathscr{B}_N^i$ and $\mathscr{B}_N^{i-1}$ overlap on half of their vertices. The blocks at stage \(n\) are then given by the blocks \(\mathscr{B}_{K_n}^i\) for $i\in\Z$. Note that all stage \(n\) blocks have the same marginal distribution by point stationarity~\cite{Thorisson2000}. { However, two neighbouring blocks are not independent due to the overlapping property. Conversely, if two blocks do not intersect, e.g.\ the blocks \(\mathscr{B}_{K_n}^i\) and \(\mathscr{B}_{K_{n}}^{i+2}\), then they, and specifically the subgraphs induced by them, are independent as a result of the renewal property of $\eta$ and the independence of the marks.}

\subsection{Connecting vertex sets that are far apart}\label{sec:LongConnect}
To make sure that the strategy outlined at the beginning of this section works and that a \(K_{n+1}\)-block at stage \(n+1\) contains a `large' connected component (we will specify this shortly), it is necessary that two stage \(n\) blocks at a given distance are connected with a sufficiently high probability to overcome potentially bad regions. 

Recall that \(K_n=(n!)^3 K^n\) for some \(K\in\N\). For $\vartheta^*\in(0,1)$ and $n\geq 2$, we { write \(v_n:=v_n(\vartheta^*)=\vartheta^*K_{n-1}\) and} define the `leftmost' and `rightmost' parts of $\mathscr{B}_{K_n}$ as
\begin{equation*}
	\begin{aligned}
		& \mathscr{V}_{\ell}^n(\vartheta^*) := \set{\X_{-K_n},\dots,\X_{-K_n+\lfloor v_n\rfloor-1}} \quad \text{ and } \\ & \mathscr{V}_{r}^n(\vartheta^*):= \set{\X_{K_n-\lfloor v_n \rfloor},\dots, \X_{K_n-1}}.
	\end{aligned}
\end{equation*}
Note that $K_{n-1}\ll K_n$ and so $\mathscr{V}_{\ell}^n({\vartheta^*})$ (resp.\ $\mathscr{V}_{r}^n(\vartheta^*)$) is only a relatively small number of vertices at the very left (resp.\ right) end of the block $\mathscr{B}_{K_n}$. Before calculating the probability of the two sets $\mathscr{V}_{\ell}^n({\vartheta^*})$ and \(\mathscr{V}_{r}^n({\vartheta^*})\) being connected, we need to understand the behaviour of the vertex marks inside each set. 

For \(\mu\in(0,\nicefrac{1}{2})\), we denote for all \(i=1,\dots, {\lfloor v_n^{1-\mu}\rfloor}\) by
\begin{equation*}
	N_\ell^n(i):=\sum_{S\colon(X,S)\in \Vl} \mathbbm{1}_{\set[\big]{S\leq \tfrac{i}{\lfloor v_n^{1-\mu}\rfloor}}}
\end{equation*}		
the empirical mark counts in \(\mathscr{V}_{\ell}^n({\vartheta^*})\). We say further that \(\mathscr{V}_{\ell}^n({\vartheta^*})\) has \emph{\(\mu\)-regular} vertex marks if 
\begin{equation*}
	\begin{aligned}
		N_\ell^n(i) \geq \frac{v_n i}{2\lfloor v_n ^{1-\mu}\rfloor}
	\end{aligned}
\end{equation*}
for all \(i=1,\dots, \lfloor v_n^{1-\mu}\rfloor\). A simple calculation yields that
\[\E N_\ell^n(i)=\frac{v_n i}{\lfloor v_n^{1-\mu}\rfloor}.\]
Hence, by the Chernoff bound for Uniform\((0,1)\) random variables, we have
\begin{equation*}
	\P\set[\big]{N_\ell^n(i)<\E[N_\ell^n(i)]/2} \leq \exp\big({-c \, i \, v_n^\mu}\big) 
\end{equation*}
for some constant \(c>0\), independent of all model parameters. Therefore,
\begin{equation}
	\P\{\Vl \text{ is } \mu\text{-regular}\}\geq 1-K_{n-1}^{1-\mu}\exp(-c v_n^\mu), \label{eqMuReg}
\end{equation}
and the same holds verbatim for \(\mathscr{V}_{r}^n(\vartheta^*)\). 

Consequently, both sets are \(\mu\)-regular with a stretched exponential error bound already in the first stage for a sufficient large \(K\). We therefore focus on the case when both sets have \(\mu\)-regular vertex marks when calculating the probability of both sets being connected, which we do in the following lemma that also illuminates the necessity of Assumption~\eqref{eqInfinite}. We take into account here that, by the evenly spaced property, the two sets are typically at distance roughly \(K_n\). { More precisely, we restrict ourselves in the following to \(K_n\)-properly spaced configurations, see~\eqref{eq:propspaced}. As the current stage \(n\), or \(K_n\) respectively, is mostly clear from context, we simply refer to the property as \emph{properly spaced}. We further fix throughout some \(a>2\lambda\) for which Definition~\ref{DefEvenlySpaced} applies.}

We denote by $\{\Vl\sim\Vr\}$ the event that the two sets are at graph distance one of each other, i.e., that there exist $\X\in\Vl$ and $\mathbf{Y}\in\Vr$ such that $\X\sim\mathbf{Y}$. We write \(\{\Vl\not\sim\Vr\}\) for the event that there is no such edge. Furthermore, we say that two vertices \((X_i,T_i)\) and \((X_j,T_j)\) are \emph{\(K_n\)-strongly connected} if their corresponding edge mark satisfies
\begin{equation}
	U_{i,j}\leq 1-\exp\big(-\rho(\beta^{-1}g(T_i,T_j)a K_n)\big), \label{eq:strConnect}
\end{equation}
We denote this event by \(\{\X_i\overset{K_n}{\rightleftharpoons}\X_j\}\) and also extend the notation to sets of vertices in the same fashion as `$\sim$'. Since \(1-e^{-x}<x\), two \(K_n\)-properly spaced vertices that are \(K_n\)-strongly connected are always connected by an edge in $\cG_\beta$. Again, if the scale \(n\) is clear from the context we simply write \(\{\X_i{\rightleftharpoons}\X_j\}\) and say that the vertices are \emph{strongly connected}.

\begin{lemma}\label{lemLongConnection}
	Let $\vartheta^*\in(0,1)$ and { recall $v_n=\vartheta^*K_{n-1}$. 
		We have for all $\mu\in(0,\nicefrac{1}{2})$} and $n\in\N$,  
	\begin{equation*}
		\begin{aligned}
			\P\Big(\Vl & \not\sim \Vr \ \big| \  \{\Vl \text{ and } \Vr \text{ are } \mu\text{-regular}\}, \{|X_{-K_n}-X_{-K_n -1}|\leq a  K_n\} \Big)  \\ 
			& \leq \P\Big(\Vl\overset{K_n}{\not\rightleftharpoons}\Vr \, \Big| \, \{\Vl \text{ and }\Vr \text{ are }\mu\text{-regular}\}\Big)
			\\
			&\leq \exp\Big(-\frac{v_n^2}{4}\int_{2 v_n^{\mu-1}}^{ 1} \d s \int_{2 v_n^{\mu-1}}^{1} \d t \ \rho\big(\beta^{-1}g(s,t)a K_n\big)\Big). 
		\end{aligned}
	\end{equation*}
\end{lemma}

\begin{proof}
	To lighten notation, we write $\mathscr{V}_\ell=\Vl$, and $\mathscr{V}_r=\Vr$, and denote 
	\[
	\cE:= \cE(n) =\{\mathscr{V}_\ell \text{ is } \mu\text{-regular}\}\cap\{\mathscr{V}_r \text{ is } \mu\text{-regular}\}. 
	\] 
	Denote by $F_\ell$ and \(F_r\) the empirical distribution function of the vertex marks corresponding to $\mathscr{V}_\ell$ and \(\mathscr{V}_r\) respectively. Writing $h:=\lfloor v_n^{1-\mu}\rfloor$, we have on the event $\cE$ for $t\in[0,1]$ by the definition of $\mu$-regularity
	\begin{equation}\label{eqEmpDF}
		\begin{aligned}
			v F_\ell(t) 
			& 
			= \sum_{i=-K_n}^{-K_n+\lfloor v_n\rfloor-1} \mathbbm{1}_{\{T_i\leq t\}} 
			\geq  N_\ell^n(\lfloor th\rfloor)\geq \frac{v_n\lfloor th\rfloor}{2h}\geq \frac{v_n}{2}(t-1/h)
			\\ &
			\geq \frac{v_n}{2}(t-2v_n^{\mu-1}). 
		\end{aligned}
	\end{equation}
	{ Let now 
		\[
		G(t) =
		\begin{cases}
			0, & t<2v_n^{\mu-1}, 
			\\
			\tfrac{1}{2}(t-2v_n^{\mu-1}), & 2v_n^{\mu-1}\leq t\leq 1,
			\\
			\tfrac{1}{2}-v_n^{\mu-1}, & t>1,
		\end{cases}
		\]
		which defines a measure supported on \([2v_n^{\mu-1},1]\) with \(\d G(t) = \nicefrac{\d t}{2}\). Using~\eqref{eqEmpDF}, we thus obtain on the event of \(\mu\)-regularity for any bounded function \(\rho\colon[0,1]\to [0,\infty)\) that
		\begin{equation}\label{eq:boundEmpSum}
			\begin{aligned}
				\sum_{i=-K_n}^{-K_n+\lfloor v_n \rfloor -1} \rho(T_i) = \int_0^1 \rho(t) \, \d F_\ell(t) \geq \int_{2v_n^{\mu-1}}^1 \tfrac{1}{2}\rho(t) \, \d t.
			\end{aligned}
		\end{equation}
		The same applies verbatim to $F_r$.
	}%

	Next, we write \(\P^{\eta_0}\) for the probability measure \(\P\) given the vertex locations \(\eta_0\) which, by construction, is a product measure with \(\operatorname{Uniform}(0,1)\) marginals. Recall that \(\mathbf{P}_0\) denotes the law of \(\eta_0\). We then have
	\begin{equation*}
		\begin{aligned}
			\P \Big( \Vl & \not\sim \Vr \, \Big| \,  \cE\cap \{|X_{-K_n}-X_{-K_n -1}|\leq a  K_n\} \Big) \\
			& = \frac{\int \1_{\{\omega\text{ properly spaced}\}} \P^{\eta_0=\omega}\big(\Vl\not\sim\Vr \,\big| \, \cE\big) \mathbf{P}_0(\d\omega)}{\mathbf{P}_0(\eta_0 \text{ properly spaced})}.
		\end{aligned}
	\end{equation*}
	We focus on the inner probability in the numerator's integral. By construction, under \(\P^{\eta_0=\omega}\), two vertices $(x_i,T_i)\in \mathscr{V}_\ell$ and $(x_j,T_j)\in \mathscr{V}_r$ are connected, whenever they are strongly connected, i.e.\ their corresponding edge mark satisfies~\eqref{eq:strConnect},
	since $\omega$ is properly spaced. Thus, the first inequality of the lemma follows immediately. In particular, there always exists an edge connecting $\mathscr{V}_\ell$ and $\mathscr{V}_r$ if
	\[
	\Sigma:=\sum_{\substack{(x_i,T_i)\in \mathscr{V}_\ell,\\ (x_j,T_j)\in \mathscr{V}_r}} \mathbbm{1}_{\big\{{U}_{i,j}\leq 1-\exp(-\rho(\beta^{-1}g(T_i,T_j)a  K_n)\big\}}>0.
	\] 
	Since the edge marks are independent of the vertex marks and locations, we have
	\begin{equation*}
		\begin{aligned}
			\E^{\eta_0=\omega}(\1_{\{\Sigma = 0\}}\1_\cE) 
			& 
			\leq \E^{\eta_0=\omega}\Bigg[\1_\cE\prod_{\substack{(x,T)\in \mathscr{V}_\ell \\ (y,S)\in \mathscr{V}_r}}\exp\big(-\rho(\beta^{-1}g(T,S)a  K_n)\big) \Bigg]  
			\\ & 
			= \E^{\eta_0=\omega}\Bigg[\exp\Big(-\sum_{\substack{(x,T)\in \mathscr{V}_\ell \\ (y,S)\in \mathscr{V}_r}} \rho\big(\beta^{-1}g(T,S)a  K_n\big)\Big) \1_\cE\Bigg] 
			\\ &
			\leq \E^{\eta_0=\omega}\Bigg[\exp\Big(-\frac{v_n^2}{4}\int_{2v_n^{\mu-1}}^{1} \d t \int_{2v_n^{\mu-1}}^{1} \d s \ \rho\left(\beta^{-1}g(t,s)a  K_n\right)\Big) \1_\cE\Bigg]
			\\ &
			= {\P}^{\eta_0=\omega}(\cE) \, \exp\Big(-\frac{v_n^2}{4}\int_{2 v_n^{\mu-1}}^{1} \d t \int_{2 v_n^{\mu-1}}^{1} \d s \ \rho\left(\beta^{-1}g(t,s)a  K_n\right)\Big),
		\end{aligned}
	\end{equation*}
	where we used~\eqref{eq:boundEmpSum} twice in the second to last step. 
	The proof is concluded by the observation that the established bound is uniform in all properly spaced configurations $\omega$.
\end{proof}

\subsection{Renormalisation scheme}
For \(N\in\N, i\in\Z\), we denote by $\C_\beta(\mathscr{B}_N^i)$ the largest connected component of the subgraph of $\cG_\beta$ on the vertices of $\mathscr{B}_N^i$. 
For $\vartheta\in(0,1)$, we say a block $\mathscr{B}_N^i$ is \emph{$\vartheta$-good}, if it contains a connected component of size at least $2\vartheta N$; otherwise we call it \emph{$\vartheta$-bad}. We denote by
\begin{equation*}
	p_\beta(N, \vartheta):= \P\big\{\sharp\C_\beta(\mathscr{B}_N)< 2N\vartheta\big\} = \P\{\mathscr{B}_N \text{ is }\vartheta\text{-bad}\}
\end{equation*}
the probability that the block $\mathscr{B}_N$ is $\vartheta$-bad. We will show that the probability of $\mathscr{B}_{K_n}$ being $\vartheta$-bad can be bounded by the probability that the smaller block $\mathscr{B}_{K_{n-1}}$ is bad with a slightly larger value of $\vartheta$. 

As a first step, we show that the subgraph of $\cG_\beta$ induced by $\mathscr{B}_{K_{n+1}}$ typically contains a connected component of volume proportion at least $\vartheta-\varepsilon$, whenever $\mathscr{B}_{K_n}$ is $\vartheta$-good with sufficiently large probability. This is an adaptation of \cite[Lemma~2]{DCGT2019} to our setting. Here, we have to deal with the positive correlations between clusters to make use of Lemma~\ref{lemLongConnection}, cf.~\eqref{eqIntuition}, and Lemma~\ref{lem:Decorrelate}. Afterwards, we show that for sufficiently large initial scales $K$, the subgraph induced by $\mathscr{B}_{K_1}$ contains a large cluster whenever \(\beta\) is large enough. Combining both results yields Proposition~\ref{thmInfinite}.

Recall that we have $K_n=(n!)^3 K^n$. Define the sequence $(C_n)_{n\in\N}$ by setting $C_n=n^3 K$. Then $K=K_1=C_1$ and $K_n=C_n K_{n-1}$ for $n\geq 2$.

\begin{lemma} \label{lemProbBadIneq} Let $\vartheta^*\in(\nicefrac{3}{4},1)$ and $\vartheta\in(\vartheta^*,1)$. Under the assumptions of Proposition~\ref{thmInfinite} there exists $M>0$ such that for all $K\geq M$ and $n\geq 2$, we have $\vartheta-\nicefrac{2}{C_n}\geq \vartheta^*$, and 
	\begin{equation*}
		p_\beta(K_n,\vartheta-\nicefrac{2}{C_n}) \leq \tfrac{1}{100}p_\beta (K_{n-1},\vartheta)+2 C_n^2p_\beta(K_{n-1},\vartheta)^2. 
	\end{equation*}	 
\end{lemma}

\begin{proof}
	Let $\vartheta':=\vartheta-\nicefrac{2}{C_n}$. Consider the blocks $\mathscr{B}_{K_{n-1}}^i$ for $|i|\in\{0,\dots, C_n-1\}$, which together form $\mathscr{B}_{K_n}$, and their largest connected components $\C_\beta(\mathscr{B}_{K_{n-1}}^i)$. Since $\vartheta>\nicefrac{3}{4}$, the cluster $\C_\beta(\mathscr{B}_{K_{n-1}}^i)$ is unique if $\mathscr{B}_{K_{n-1}}^i$ is $\vartheta$-good.  Furthermore, due to the overlapping property of neighbouring blocks, the largest components of two adjacent $\vartheta$-good blocks have to intersect in at least one vertex. Hence, if all the blocks $\mathscr{B}_{K_{n-1}}^i$ are $\vartheta$-good, then $\mathscr{B}_{K_{n}}$ is $\vartheta$-good as well.
	
	Define now for every $|i|\in\{0,\dots, C_n-1\}$ the event
	\[
	\cE_i:=\cE_i(n,\beta) = \big\{\sharp\C_\beta(\mathscr{B}_{K_{n-1}}^i)<2\vartheta K_{n-1} \big\} \cap \Bigg[\bigcap_{\substack{|j|=0 \\ j\not\in\{i-1,i,i+1\}}}^{C_n-1} \big\{\sharp\C_\beta(\mathscr{B}_{K_{n-1}}^j)\geq 2\vartheta K_{n-1}\big\}\Bigg].
	\]
	That is, the block $\mathscr{B}_{K_{n-1}}^i$ is $\vartheta$-bad but all blocks $\mathscr{B}_{K_{n-1}}^j$ which it does not intersect are $\vartheta$-good. If we write
	\[
	\cF_i:=\cF_i(n,\beta) := \cE_i \cap \{\mathscr{B}_{K_n} \text{ is }\vartheta{{'}}\text{-bad}\},
	\]
	then $\mathscr{B}_{K_n}$ being $\vartheta'$-bad implies that either $\cF_i$ occurs for some $i$ or at least two disjoint stage \((n-1)\) blocks are $\vartheta$-bad, since otherwise $\mathscr{B}_{K_n}$ would be $\vartheta$- and therefore also $\vartheta'$-good. Consequently,
	\begin{equation*}
		\begin{aligned}
			p_\beta(K_n, \vartheta') 
			& 
			= p_\beta(C_nK_{n-1}, \vartheta')  \leq \sum_{|i|=0}^{C_n-1} {\P}(\cF_i) + \binom{C_n}{2} p_{\beta}(K_{n-1}, \vartheta)^2 
			\\ & 
			\leq p_\beta(K_{n-1},\vartheta) \sum_{|i|=0}^{C_n-1} {\P}(B_{C_nK_{n-1}} \text{ is } \vartheta'\text{-bad}\mid\cE_i) +2C_n^2 p_\beta(K_{n-1},\vartheta)^2,
		\end{aligned}
	\end{equation*}
	{ using { in the first inequality} the independence of subgraphs on disjoint blocks.}
	To finish the proof it therefore remains to bound the sum of the conditional probabilities by $\nicefrac{1}{100}$. To this end, define 
	\begin{equation}\label{eq:largestClust}
		\mathscr{C}_i^\ell := \C_{i}^\ell(n,\beta)=\bigcup_{j=1-C_n}^{i-2} \C_\beta(\mathscr{B}_{K_{n-1}}^j) \ \text{ and } \ \C_i^r := \C_i^r(n,\beta) = \bigcup_{j=i+2}^{C_n-1} \C_\beta(\mathscr{B}_{K_{n-1}}^j),
	\end{equation}
	{ the union of all largest clusters to the left and to the right, respectively, of the bad block \(\mathscr{B}_{K_{n-1}}^i\).} 
	Conditioned on $\cE_i$, both sets $\C_i^\ell$ and $\C_i^r$ are connected sets. Furthermore, if $i\in\{C_n-2, C_n-1\}$, then
	\begin{equation*}
		\begin{aligned}
			\sharp\C_i^\ell \geq 2(C_n-2)\vartheta K_{n-1} \geq 2\vartheta' C_n K_{n-1}
		\end{aligned}
	\end{equation*}
	and hence $\mathscr{B}_{K_n}$ is $\vartheta'$-good. 
	\begin{figure}
		\begin{center}
			\begin{tikzpicture}[scale = 0.8,every node/.style={scale=0.7},arrow/.style={-latex, shorten >=1ex, shorten <=1ex, bend angle=45}]
				\draw (0,0) rectangle (3,2)[thick]{};
				\draw (1.5,0) rectangle (4.5,2)[]{};
				\draw (3,0) rectangle (6,2)[thick]{};
				\draw (4.5,0) rectangle (7.5,2)[fill = lightgray]{};
				\draw (7.5,0) rectangle (10.5,2)[fill = lightgray]{};
				\draw (6,0) rectangle (9,2)[thick, fill = gray]{};
				
				\draw (9,0) rectangle (12,2)[thick]{};
				\draw (10.5,0) rectangle (13.5,2)[]{};
				
				
				\draw[decorate, decoration = {snake, segment length = 15 pt, amplitude = 5mm}, thick, color= blue] (0.5, 1)--(4.5,1);
				
				\draw[decorate, decoration = {snake, segment length = 15 pt, amplitude = 5mm}, thick, color= red] (13, 1)--(10.5,1);
				
				\draw[bend angle = 60, bend left, dashed, thick] (2.5,1.5) to (11,1.5);
			\end{tikzpicture}
			\caption{The overlapping blocks of scale \(n-1\) that together form the scale \(n\) block. The cluster \(\C^-_i\) (resp.\ \(\C^+_i\)) on the left in blue (resp.\ on the right in red). The dark block is the bad block and in light gray are the non overlapping halves of the two neighbouring blocks. The dotted line indicates the existence of an edge connecting \(\C^-_i\) and \(\C^+_i\) avoiding the bad region.} 
			\label{fig:SketchofProof}
		\end{center}
	\end{figure}
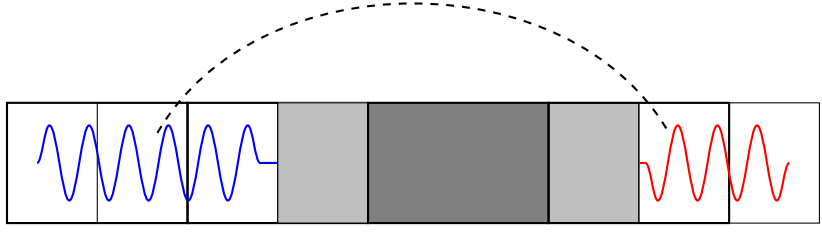
	The same holds true for $\C_i^r$ if $i\in\{1-C_n, 2-C_n\}$. Therefore, the bad block and any neighbouring block cannot be the left- or the right-most one in $\mathscr{B}_{K_n}$. This then guarantees that $\C_i^\ell,\C_i^r\neq\emptyset$. Further, if $\C_i^\ell$ and $\C_i^r$ are connected directly by an edge, we have
	\[\sharp\C_i^\ell +\sharp\C_i^r\geq \vartheta K_{n-1}(C_n+i-2)+\vartheta K_{n-1}(C_n-i-2)=2\vartheta K_{n-1}(C_n-2)\geq 2\vartheta' C_n K_{n-1},\]
	and $\mathscr{B}_{K_n}$ is again $\vartheta'$-good, see Figure~\ref{fig:SketchofProof}; {thus,
		\begin{equation*}
			\begin{aligned}
				{\P}(\mathscr{B}_{K_{n}} \text{ is }\vartheta'\text{-bad}\mid\cE_i) & \leq {\P}(\C^\ell_i\not\sim\C^r_i\mid\cE_i).
			\end{aligned}
		\end{equation*}    
		
		To control the probability on the right-hand side, we make now use of the evenly spaced property. Recall that on a properly spaced configuration a strong connection in the sense of~\eqref{eq:strConnect} is harder to achieve than a normal connection, which yields
		\begin{equation}\label{eq:propSpacedBound}
			\begin{aligned}
				\P(\C_i^\ell \not\sim\C_i^r\mid\cE_i)
				&
				\leq \int\limits_{\{\omega \phantom{x} \substack{\text{\tiny properly}\\ \text{\tiny spaced}}\}} \P^{\eta_0=\omega}(\C_i^{\ell}\not\rightleftharpoons\C_i^r\mid \cE_i^{\omega}) \mathbf{P}_0(\d \omega)+ o(K_n^{-2}),
			\end{aligned}
		\end{equation}
		by Definition~\ref{DefEvenlySpaced}, where \(\cE_i^{\omega}\) denotes the event that \(\cE_i\) occurs on the properly spaced configuration \(\eta_0=\omega\); note that this has positive probability. Again, we have written \(\P^{\eta_0=\omega}\) for the conditional probability given the locations \(\eta_0=\omega\). Let us next define 
		\[
		\cA_i:=\cA_i(n,\beta) = \Bigg[\bigcap_{\substack{|j|=0 \\ j\not\in\{i-1,i,i+1\}}}^{C_n-1} \big\{\sharp\C_\beta(\mathscr{B}_{K_{n-1}}^j)\geq 2\vartheta K_{n-1}\big\}\Bigg],
		\]
		so that \(\cE_i=\{\sharp\C_\beta(\mathscr{B}_{K_{n-1}}^i)<2\vartheta K_{n-1} \big\} \cap\cA_i\) and denote \(\cA_i^\omega\) for the event that \(\cA_i\) occurs on the location configuration \(\eta_0=\omega\). As strong connections only depend on the involved vertex and edge marks but \emph{not} on the precise differences of the locations, we obtain for any properly spaced configuration \(\eta_0=\omega\)
		\begin{equation*}
			\begin{aligned}
				\P^{\eta_0=\omega}(\C_i^{\ell}\not\rightleftharpoons\C_i^r\mid \cE_i^{\omega}) = \P^{\eta_0=\omega}(\C_i^{\ell}\not\rightleftharpoons\C_i^r\mid \cA_i^{\omega}),
			\end{aligned}
		\end{equation*}
		as \(\C_i^\ell\) and \(\C_i^r\) do not intersect the bad block \(\mathscr{B}^i_{K_{n-1}}\) or its adjacent blocks. 
		On \(\cA_i^\omega\), the clusters \(\C_i^\ell\) and \(\C_i^r\) are the unique largest clusters on the left or on the right side of the bad block respectively, which provides positive information on the connectivity of the involved vertices. It thus seems plausible that
		\begin{equation} \label{eqIntuition}
			\begin{aligned}
				\P^{\eta_0=\omega}(\C_i^{\ell}\not\rightleftharpoons\C_i^r\mid \cA_i^{\omega}) \leq \P^{\eta_0=\omega}(\Vl\not\rightleftharpoons\Vr),
			\end{aligned}
		\end{equation}
		for which Lemma~\ref{lem:Decorrelate} below provides rigorous justification.} Let us for the moment assume that~\eqref{eqIntuition} holds true. We then infer by combing this with~\eqref{eq:propSpacedBound}, the help of Lemma~\ref{lemLongConnection}, and \(\mu\)-regularity~\eqref{eqMuReg}
	\begin{equation*}
		\begin{aligned}
			\sum_{|i|=0}^{C_n-1} & {\P}(B_{C_nK_{n-1}} \text{ is }\vartheta'\text{-bad}\mid\cE_i) \\ 
			& \leq 2n^3 C_1\Big[ \exp\Big(-Cv^2\int_{[v^{\mu-1}, 1-v^{\mu-1}]^2}\d(t,s) \ \rho\big(\beta^{-1}g(s,t)aK_n\big)\Big) \\ 
			& \phantom{02n^3C_1\Big[\exp }+ K_{n-1}^{1-\mu}\exp(-cK_{n-1}^\mu) + o(K_n^{-2})\Big] \\ 
			& \leq \frac{1}{100}.
		\end{aligned}
	\end{equation*}
	by Assumption~\eqref{eqInfinite} for sufficiently large \(K=C_1\), as desired.
\end{proof}

It remains to justify~\eqref{eqIntuition} to formally conclude the proof of Lemma~\ref{lemProbBadIneq}.

\begin{lemma}\label{lem:Decorrelate} 
	For all \(\beta>0, n\geq 2\), \(|i|\in\{0,\dots,C_n-1\}\) and \(K_n\)-properly spaced location configuration \(\eta_0=\omega\) we have 
	\begin{equation*}
		\begin{aligned}
			\P^{\eta_0=\omega}(\C_i^{\ell}\not\rightleftharpoons\C_i^r\mid \cA_i^{\omega}(n,\beta)) \leq \P^{\eta_0=\omega}(\Vl\not\rightleftharpoons\Vr).
		\end{aligned}
	\end{equation*}	 
\end{lemma}

\begin{proof}
	{ Recall that \(\lfloor v_n\rfloor =\sharp\Vl=\sharp\Vr\). To keep the notation readable, we omit writing the integer part function in the following and treat \(v_n\) as a natural number.} The idea is to bound the probability on the left by first uniformly choosing subsets of size \(v_n\) among the vertices of \(\C_i^\ell\) and \(\C_i^r\) and only check whether these subsets are strongly connected or not and then compare this with the case when $v_n$ vertices are chosen uniformly among \emph{all} vertices on the left and right of the box \(\mathscr{B}_{K_{n-1}}^i\).  To do so rigorously, we first have to extend our probability space. Let 
	\[
	\mathcal{L}_i:=\mathcal{L}_i(n)=\{\ell_1,\dots, \ell_{v_n}\}\subset \{-C_n K_{n-1},\dots, K_{n-1}(i-1)-1\}
	\]	
	be indices chosen among all indices of the vertices on the left side of the block \(\mathscr{B}_{K_{n-1}}^i\) uniformly without replacement, independently of everything else, and ordered from smallest absolute value to largest.  Similarly, let
	\[
	\mathcal{R}_i:=\mathcal{R}_i(n)=\{r_1,\dots,r_{v_n}\}\subset \{K_{n-1}(i+1),\dots, C_nK_{n-1}-1\}
	\]
	be another independent set of indices chosen uniformly among the indices of the vertices on the right side. We call each index in \(\cL_i\) and \(\cR_i\) \emph{tagged}. Note that the sets of indices we are sampling from are deterministic. Furthermore, we can define a joint probability measure \(\overline{\P}^{\eta_0=\omega}\) for \((\xi_0,(\cL_i,\cR_i))\), given \(\eta=\omega\) with marginal distribution \(\P^{\eta_0=\omega}\), when integrated with respect to the tags \(\cL_i\) and \(\cR_i\). 
	%
	
	Recall the definitions of \(\C_i^\ell\) and \(\C_i^r\) as the union of the largest clusters in the blocks left and right of \(\mathscr{B}_{K_{n-1}}^i\), cf.~\eqref{eq:largestClust}. To make this definition unique, we define from here on onwards, \(\C_\beta(\mathscr{B}_{K_{n-1}}^j)\) to be the largest cluster in \(\mathscr{B}_{K_{n-1}}^j\) that contains the vertex with smallest mark if the largest cluster is not unique. Then, \(\C_i^\ell\) and \(\C_i^r\) are always uniquely determined (almost surely) but on \(\cA_i^\omega\) nothing has changed. We further denote by
	\[
	\mathcal{I}_i^{\ell}:=\mathcal{I}_i^{\ell}(n,\beta)=\{j\in\Z\colon \X_j\in\C_i^\ell\} \ \text{ and } \ \mathcal{I}_i^{r}:=\mathcal{I}_i^{r}(n,\beta)=\{j\in\Z\colon \X_j\in\C_i^r\},
	\]
	the set of indices belonging to \(\C_i^\ell\) and \(\C_i^r\) respectively; we consider both sets to be ordered from smallest absolute value to largest. 
	Let us further introduce another \emph{configuration dependent} tagging. Let \(\mathbf{L}_i:=\mathbf{L}_i(n)\) be a set of \(v_n\)-many indices chosen uniformly without replacement from \(\mathcal{I}_i^{\ell}\), ordered from smallest absolute value to largest, and \(\mathbf{R}_i:=\mathbf{R}_i(n)\) the same but chosen among the indices in \(\mathcal{I}_i^{r}\), independently of \(\mathbf{L}_i\). If \(\sharp\mathcal{I}_i^{\ell}<v_n\) (resp.\ \(\sharp\mathcal{I}_i^{r}<v_n\)), we simply set \(\mathbf{L}_i=\emptyset\) (resp.\ \(\mathbf{R}_i=\emptyset\)). Let us denote the joint distribution, given \(\eta_0=\omega\), of \(\xi_0\) and \((\mathbf{L}_i,\mathbf{R}_i)\) by \(\widetilde{\P}^{\eta_0=\omega}\). Observe that on \(\cA_i^\omega\), both \(\mathcal{I}_i^{\ell}\) and \(\mathcal{I}_i^{r}\) always contain at least \(v_n\)-many elements.
	Now, by the independence of \((\cL_i,\cR_i)\) from the vertex and edge marks and by the fact that uniformly sampling from a finite set \(S\) conditioned to be contained in an independently generated subset \(S'\subset S\) has the same distribution as uniformly sampling from \(S'\), we have
	\begin{equation}\label{eq:condtionalSampling}
        \begin{aligned}
           		\widetilde{\P}^{\eta_0=\omega} 
                &
                    \big((\xi_0,(\mathbf{L}_i,\mathbf{R}_i))\in\cdot \, \big| \,\cA_i^\omega\big) 
                \\ &
                    = \overline{\P}^{\eta_0=\omega}\big((\xi_0,(\cL_i,\cR_i))\in\cdot \, \big| \,\cA_i^\omega,\{\cL_i\subset\mathcal{I}_i^{\ell}\},\{\mathcal{R}_i\subset\mathcal{I}_i^{+}\}\big). 
        \end{aligned}
	\end{equation}
	This then implies
	\begin{equation}\label{eq:fkgCalc1}
		\begin{aligned}
			\widetilde{\P}^{\eta_0=\omega}(\C_i^\ell\not\rightleftharpoons\C_i^r\mid\cA_i^\omega)  
			& 
			\leq \widetilde{\P}^{\eta_0=\omega}\Big(\bigcap_{\substack{\ell\in\mathbf{L}_i, \,  r\in\mathbf{R}_i}}\{\X_\ell\not\rightleftharpoons\X_r\} \, \Big| \,\cA_i^\omega\Big) 
			\\ & 
			=\overline{\P}^{\eta_0=\omega}\Big(\bigcap_{\ell\in\cL_i, \, r\in\cR_i}\{\X_\ell\not\rightleftharpoons\X_r\} \, \Big| \,\cA_i^\omega,\{\cL_i\subset\mathcal{I}_i^{\ell}\},\{\cR_i\subset\mathcal{I}_i^{r}\}\Big) 
			\\ & 
			= \frac{\overline{\P}^{\eta_0=\omega}\Big(\bigcap\limits_{\ell\in\cL_i, \, r\in\cR_i}\{\X_\ell\not\rightleftharpoons\X_r\}, \cA_i^\omega, \{\cL_i\subset\mathcal{I}_i^{\ell}\}, \{\cR_i\subset\mathcal{I}_i^{r}\}\Big) }{\overline{\P}^{\eta_0=\omega}(A_i^\omega, \{\cL_i\subset\mathcal{I}_i^{\ell}\},\{\cR_i\subset\mathcal{I}_i^{r}\})}.
		\end{aligned}
	\end{equation}
	The probability in the last line's numerator can be written using the independence of \((\cL_i,\cR_i)\) from the graph as
	\begin{equation}
		\begin{aligned}\label{eq:fkgCalc2}
			& \overline{\P}^{\eta_0=\omega}  \Big(\bigcap\limits_{\ell\in\cL_i, \, r\in\cR_i}\{\X_\ell\not\rightleftharpoons\X_r\}, \cA_i^\omega, \{\cL_i\subset\mathcal{I}_i^{\ell}\}, \{\cR_i\subset\mathcal{I}_i^{r}\}\Big) \\ & \quad  = \sum_{L,R} \P^{\eta_0=\omega}\Big(\bigcap\limits_{\ell\in L, \, r\in R}\{\X_\ell\not\rightleftharpoons\X_r\}, \cA_i^\omega, \{L\subset\mathcal{I}_i^{\ell}\}, \{R\subset\mathcal{I}_i^{r}\}\Big)P\{\cL_i=L, \cR_i=R\},
		\end{aligned}
	\end{equation}
	where the summation runs over all subsets of size \(v_n\) of the sets, from which \(\cL_i\) and \(\cR_i\) are drawn, and we have written \(P\) for the law of \((\cL_i,\cR_i)\).
	Now, \(\P^{\eta_0=\omega}\) is by construction a product measure with \(\operatorname{Uniform}(0,1)\) marginals (the vertex and edge marks). Moreover, the event \(\bigcap_{\ell,r}\{\X_\ell\not\rightleftharpoons\X_r\}\) is clearly decreasing whenever involved vertex and/or edge marks are decreased (and potentially new edges are added to the graph). On the contrary, the event \(A_i^\omega\cap\{L\subset\mathcal{I}_i^{\ell}\}\cap\{R\subset\mathcal{I}_i^{r}\}\) is increasing whenever vertex and/or edge marks are decreased as additional edges only increase the size of the largest cluster. Hence, good boxes always remain good when an edge is added. Also, since, on \(\cA_i^\omega\), the largest clusters on the left and right are uniquely determined and contain at least a \(\vartheta>\nicefrac{3}{4}\) proportion of the vertices, an additional edge can only make the sets \(\mathcal{I}_i^{\ell},\mathcal{I}_i^{r}\) larger and whenever \(\{L\subset\mathcal{I}_i^{\ell}\}\cap\{R\subset\mathcal{I}_i^{r}\}\) occurs without the newly added edge, it also occurs after the edge has been added. Hence, we can apply the FKG-inequality~\cite{FKG71} and infer
	\begin{equation*}
		\begin{aligned}
			\P^{\eta_0=\omega}\Big(&\bigcap\limits_{\ell\in L, \, r\in R}\{\X_\ell\not\rightleftharpoons\X_r\}, \cA_i^\omega, \{L\subset\mathcal{I}_i^{\ell}\}, \{R\subset\mathcal{I}_i^{r}\}\Big) \\ &\leq \P^{\eta_0=\omega}\Big(\bigcap\limits_{\ell\in L, \, r\in R}\{\X_\ell\not\rightleftharpoons\X_r\}\Big) \P^{\eta_0=\omega}\big(A_i^\omega, \{L\subset\mathcal{I}_i^{\ell}\}, \{R\subset\mathcal{I}_i^{r}\}\big) \\
			& =\P^{\eta_0=\omega}\Big(\bigcap\limits_{\ell=-K_n}^{-K_n+v-1} \bigcap_{ r=K_n-v}^{K_n-1}\{\X_\ell\not\rightleftharpoons\X_r\}\Big) \P^{\eta_0=\omega}\big(A_i^\omega, \{L\subset\mathcal{I}_i^{\ell}\}, \{R\subset\mathcal{I}_i^{r}\}\big),
		\end{aligned}
	\end{equation*}
	where the last equality follows from the fact that under \(\P^{\eta_0=\omega}\) all vertex and edge marks are i.i.d.\ and the event of two vertices being strongly connected is independent of the precise spatial location, cf.\ \eqref{eq:strConnect}. Plugging this back into~\eqref{eq:fkgCalc2}, we infer
	\begin{align*}
			& 
			\overline{\P}^{\eta_0=\omega}  \Big(\bigcap\limits_{\ell\in\cL_i, \, r\in\cR_i}\{\X_\ell\not\rightleftharpoons\X_r\}, \cA_i^\omega, \{\cL_i\subset\mathcal{I}_i^{\ell}\}, \{\cR_i\subset\mathcal{I}_i^{r}\}\Big) 
			\\ & 
			\quad  \leq \P^{\eta_0=\omega}\Big(\bigcap\limits_{\ell=-K_n}^{-K_n+v-1} \bigcap_{ r=K_n-v}^{K_n-1}\{\X_\ell\not\rightleftharpoons\X_r\}\Big) 
			\\ & 
			\phantom{\P^{\eta_0=\omega}(\bigcap}\times \sum_{L,R}\P^{\eta_0=\omega}\big(\cA_i^\omega, \{L\subset\mathcal{I}_i^{\ell}\}, \{R\subset\mathcal{I}_i^{r}\}\big)P\{\cL_i=L,\cR_i=R\} \\
			&\quad = \P^{\eta_0=\omega}\Big(\bigcap\limits_{\ell=-K_n}^{-K_n+v-1} \bigcap_{ r=K_n-v}^{K_n-1}\{\X_\ell\not\rightleftharpoons\X_r\}\Big) \overline{\P}^{\eta_0=\omega}\big(\cA_i^\omega, \{\cL_i\subset\mathcal{I}_i^{\ell}\},\{\cR_i\subset\mathcal{I}_i^{r}\}\big).
	\end{align*}
	Combined with~\eqref{eq:fkgCalc1}, this finally yields
	\begin{equation*}
		\P^{\eta_0=\omega}(\C_i^\ell\not\rightleftharpoons\C_i^r\mid\cA_i^\omega) \leq \P^{\eta_0=\omega}\{\Vl\not\rightleftharpoons\Vr\},
	\end{equation*}
	as desired.
\end{proof}

We close the section with the following lemma which establishes the probability bounds necessary to initialise the renormalisation scheme.

\begin{lemma}  \label{lemBox1}
	Let \(\eta\) be an evenly spaced renewal process. Then, for every kernel \(g\), every profile function \(\rho\) and every $\vartheta\in(0,1)$, there exist constants $M>0$, and $B>0$ such that for all $K=C_1>M$ and $\beta>B \, C_1$, we have
	\begin{equation*}
		p_\beta(\vartheta, C_1) \leq \tfrac{1}{400} C_1^{-2}. 
	\end{equation*}
\end{lemma}

\begin{proof}
	Denote by $\ER_n^\lambda$ an Erd\H{o}s--R\'{e}nyi-graph on $n$ vertices with edge probability $\nicefrac{\lambda}{n}$; denote its law by $\operatorname{P}_n^\lambda$. If $\lambda>1$, then $\ER_n^\lambda$ is supercritical, i.e., for all $\varepsilon_1>0$, there exists $c>0$ and $N(\varepsilon_1,\lambda)>0$ such that 
	\begin{equation}
		\operatorname{P}_n^{\lambda}\big\{\sharp\C(\ER_n^\lambda>cn\big\}\geq 1-\varepsilon_1, \quad n\geq N(\varepsilon_1,\lambda), \label{eqER}
	\end{equation}
	where $\C(\ER_n^\lambda)$ denotes the largest connected component of the graph \(\ER_n^\lambda\) \cite{vdH2017}. 
	We aim to compare this behaviour with the behaviour of the finite graph induced by the finite block \(\mathscr{B}_{C_1}\) by making use of the evenly spaced property. By Definition~\ref{DefEvenlySpaced}, we have for large enough \(C_1\), 
	\[
	p_\beta(\vartheta, C_1)\leq \P\big(\mathscr{B}_{C_1} \text{ is }\vartheta\text{-bad}\, \big| \, |X_{-C_1}-X_{C_1 -1}|\leq a C_1\big) + o(C_1^{-2}).
	\]
	In the remaining proof we hence work conditionally on the event \(|X_{-C_1}-X_{C_1 -1}|\leq a C_1\) and denote again by \(\P^{\eta_0=\omega}\) the probability measure, given a properly spaced point configuration \(\omega\), cf.~\eqref{eq:propspaced}. We assume without loss of generality that \(\rho(1)>0\). We further assume that \(g\) is bounded and remark on the unbounded case below. Now, fix \(\beta>a ||g||_\infty C_1\). Then, for all $(x_i,T_i),(x_j,T_j)\in \mathscr{B}_{C_1}$, we have
	\[\P^{\eta_0=\omega}\big\{(x_i,T_i)\sim(x_j,T_j)\big\}\geq \rho(1)\]
	and we focus on the subgraph on \(\mathscr{B}_{C_1}\) where only the edges with marks smaller than \(\rho(1)\) are present, which is now independent of vertex marks and locations.
	For a fixed $\lambda>1$, we set $c$ accordingly to above, fix $\varepsilon_2<\vartheta/c$ and choose $C_1$ large enough such that
	\[
	2\varepsilon_2 C_1\rho(1) \geq \lambda \text{ and }\lfloor 2\varepsilon_2 C_1\rfloor \geq N(\varepsilon_1, \lambda).
	\]
	Denote by $\mathscr{H}$ the subgraph on the vertices $\{\X_0,\dots, \X_{\lfloor 2\varepsilon_2 C_1\rfloor}\}\subset \mathscr{B}_{C_1}$. By \eqref{eqER}, we have
	\begin{equation*}
		\begin{aligned}
			\P^{\eta_0=\omega}\big\{\sharp\C_\beta(\mathscr{H})> c\cdot 2\varepsilon_2 C_1\big\} & \geq \operatorname{P}_{2\varepsilon_2 C_1}^{\lambda}\Big\{\sharp\C\big(\ER_{2\varepsilon_2 C_1}^\lambda\big)>c\cdot 2\varepsilon_2 C_1\Big\} \\ & \geq 1-\varepsilon_1.
		\end{aligned}
	\end{equation*}
	On $\{\sharp\C_\beta(\mathscr{H})>c\cdot 2\varepsilon_2 C_1\}$, the block $\mathscr{B}_{C_1}$ is $\vartheta$-good if enough of the remaining vertices in $\mathscr{B}_{C_1}\setminus \mathscr{H}$ are connected to $\C_\beta(\mathscr{H})$. Each such remaining vertex is connected to $\C_\beta(H)$ with a probability of at least 
	\[
	q:=q(C_1)=1-(1-\rho(1))^{2\varepsilon_2 C_1}.
	\]
	For $\psi> (\vartheta-c\varepsilon_2)/(1-\varepsilon_2)$ and $C_1$ large enough such that $q>\psi$, we have, by writing $F_{\operatorname{Bin}(n, p)}$ for the distribution function of a binomial random variable with parameters $n$ and $p$, that
	\begin{equation*}
		\begin{aligned}
			\P^{\eta_0=\omega}\big\{\sharp\C_\beta(\mathscr{B}_{C_1})>2\vartheta C_1\big\} 
			& 
			\geq (1-\varepsilon_1) \P^{\eta_0=\omega}\Big(\sharp\C_\beta(\mathscr{B}_{C_1})>2\vartheta C_1 \ \Big| \ \sharp\C_\beta(\mathscr{H})>c(2\varepsilon_2 C_1)\Big) 
			\\ & 
			\geq (1-\varepsilon_1)\Big(1-F_{\operatorname{Bin}(2(1-\varepsilon_2) C_1, q)}\big(2\psi (1-\varepsilon_2) C_1\big)\Big)
			\\ & 
			\geq 1- \exp(-c' C_1),
		\end{aligned}
	\end{equation*}
	for some $c'>0$ by a standard Chernoff bound. We conclude again with the observation that the established bound is uniform in all properly spaced configurations \(\omega\).
	
	If \(g\) is not bounded, we can argue as follows: fix a small \(\varepsilon>0\) and only consider vertices with marks smaller than \(1-\varepsilon\) and therefore each vertex is removed independently with probability \(\varepsilon\) due to independence of marks and locations. However, the new block \(\mathscr{B}_{C_1}\) still consists of order \((1-\varepsilon)C_1\) vertices with an error term exponentially small in \(C_1\) by Chernoff's bound. Furthermore, the thinned process \(\eta\) is still evenly spaced and we can repeat the proof above since it holds that \(g(s,t)\leq g(1-\varepsilon,1-\varepsilon)<\infty\) for all remaining marks \(s\) and \(t\).
\end{proof}

\subsection{Finalising the proof of Proposition~\ref{thmInfinite}} 
We are now ready to prove Proposition~\ref{thmInfinite}, which we do by following the arguments of the proof of Theorem~1(i) of \cite{DCGT2019} in the following lemma.

\begin{lemma} \label{lemProofThm}
	Let the assumptions of Proposition~\ref{thmInfinite} be fulfilled. Then there exist $\beta_c\in(0,\infty)$ such that
	\[
	\theta(\beta)={\P}\{\0\leftrightarrow\infty\text{ in }\cG_\beta\}\geq \tfrac{3}{8}
	\]
	for all $\beta>\beta_c$.
\end{lemma}

\begin{proof}
	Fix $\vartheta^*\in(\nicefrac{3}{4},1)$ and $\vartheta \in(\vartheta^*,1)$. Choose $K=C_1$ and afterwards $\beta$ both large enough, such that the assumptions of the Lemmas~\ref{lemProbBadIneq} and \ref{lemBox1} are satisfied. Recall also that $C_n=n^3 K$, $K_1=C_1$ and $K_n=C_n K_{n-1}$. Define $\vartheta_n:= \vartheta -2/C_{n+1}$ for $n\geq 2$. Since the assumptions of Lemma~\ref{lemProbBadIneq} are satisfied, we have $\vartheta_n>\vartheta^*$ for all $n$. We have by Lemma~\ref{lemBox1} that $p_\beta(C_1,\vartheta )\leq (400 C_1^2)^{-1}$, and by Lemma~\ref{lemProbBadIneq} that 
	\[
	p_\beta(K_n,\vartheta_n) \leq \tfrac{1}{100}p_\beta (K_{n-1},\vartheta_{n-1})+2 C_n^2p_\beta(K_{n-1},\vartheta_{n-1})^2, \quad \forall n\geq 2.
	\]
	Inductively, this yields $p_\beta(K_n,\vartheta_n)\leq (400 C_n^2)^{-1}$ for all \(n\in\N\), and hence
	\[{\P}\{\mathscr{B}_{K_n} \text{ is }\vartheta_n\text{-good}\} \geq 1-\tfrac{1}{400 C_n^2}\geq \tfrac{1}{2}.\]
	We derive from this
	\begin{equation*}
		\begin{aligned}
			\tfrac{3}{4} K_n & \leq \tfrac{3}{4}(2K_n) {\P}\{\mathscr{B}_{K_n} \text{ is }\nicefrac{3}{4} \text{-good}\} \leq {\E}\big[\sharp\C_\beta(\mathscr{B}_{K_n})\mathbbm{1}_{\{\mathscr{B}_{K_n} \text{ is } \nicefrac{3}{4}\text{-good}\}}\big] \\ &\leq 2 K_n {\P}_{\beta}\{\exists \text{ a cluster of size at least } \tfrac{3}{2} K_n\}.
		\end{aligned}
	\end{equation*}
	
	Dividing both sides by \(2K_n\) and then sending \(n\to\infty\) together with the translation invariance yield the desired result.
\end{proof}

In the remainder of the section we prove Corollary~\ref{cor:finiteGraphs} and specifically the existence of a component of linear size in the constructed graph sequence for large enough \(\beta\) if the assumption of Theorem~\ref{thm:main_thm}(a) is fulfilled.

\begin{proof}[Proof of Corollary \ref{cor:finiteGraphs}] 
	We are interested in the limiting behaviour of the non-negative random variable \({\sharp\mathscr{C}(\G_{n}(\beta))}/{n}\) for a fixed \(\beta\). First note that this is a translation invariant functional of the ergodic point process \(\xi\).
	Hence, \({\sharp\mathscr{C}(\G_{n}(\beta))}/{n}\) converges almost surely towards a non-negative constant.  
	Consider now for some \(K\in\N\) the subsequence of graphs \((\G_{4K_n}:n\in\N)\). Observe, that the vertices of the block \(\mathscr{B}_{K_n}\) are contained in the interval \((-2K_n, 2K_n)\) with an error term going to zero when \(K_n\to\infty\). Hence, we choose \(K\) and \(\beta\) large enough to fulfil the assumptions of the Lemmas~\ref{lemProbBadIneq} and~\ref{lemBox1} and \(n_0\) large enough such that \(\mathbf{P}_0(\mathscr{B}_{K_{n_0}}\not\subset(-2K_n,2K_n))\leq \nicefrac{1}{3}\) for all \(n\geq n_0\) and infer from the calculations in the proof of Lemma~\ref{lemProofThm}
	\[
	\P\big\{\tfrac{\sharp\mathscr{C}(\G_{4K_n}(\beta))}{4K_n}\geq \tfrac{3}{8}\big\}\geq \tfrac{1}{4}
	\]
	uniformly for \(n\geq n_0\). Hence, the considered sequence \(\big({\sharp\mathscr{C}(\G_\beta^{(n)})}/{n}\colon n\in\N\big)\) contains a subsequence with a strictly positive limit, finalising the proof.
\end{proof}

\section{Proof of {Proposition~\ref{thmFiniteComponents}}} \label{secFiniteComponent}
We define the disjoint sets of vertices
\begin{equation*}
	\begin{aligned}
		& \Gamma_k^\ell := \set{\X_{-2^k},\dots,\X_{-1}}, \qquad \Gamma_k^{\ell \ell}:= \set{\X_{-2^{k+1}},\dots, \X_{-2^{k}-1}} \\
		& \Gamma_k^r := \set{\X_0,\dots, \X_{2^k-1}}, \qquad \Gamma_k^{r r}:=\set{\X_{2^k},\dots, \X_{2^{k+1}-1}}
	\end{aligned}
\end{equation*}  
for each \(k\in\N\). We say that a \emph{crossing of the origin} occurs at stage
\begin{description}
	\item[\(k=1,\)]  if any edge connects the sets \(\Gamma ^{\ell}\cup\Gamma ^{\ell\ell}\) and \(\Gamma ^{r}\cup\Gamma ^{rr}\) or at stage
	\item[\(k\geq 2,\)] if any edge connects either \(\Gamma_k^{\ell\ell}\) to \(\Gamma_k^{rr}\), \(\Gamma_k^{\ell\ell}\) to \(\Gamma_k^r\) or \(\Gamma_k^\ell\) to \(\Gamma_k^{rr}\). Note that any edges between \(\Gamma_k^\ell\) and \(\Gamma_k^r\) have by necessity already been considered at an earlier stage.
\end{description}
We denote by \(\chi(k)\in\{0,1\}\) the indicator of the event that a crossing of the origin occurs at stage \(k\in\N\). The event that there is no edge crossing the origin is then given by \(\bigcap_k\{\chi(k)=0\}\). We write
\begin{equation} \label{eq:noCrossing}
	\P\Big\{\bigcap_{k\in\N}\chi(k)=0\Big\} = \mathbf{E}_0\Big[\P^{\eta_0}\Big\{\bigcap_{k\in\N}\chi(k)=0\Big\}\Big],
\end{equation}
an focus on the right-hand side's inner probability where again \(\mathbf{E}_0\) denotes the expectation of \(\eta_0\) and \(\P^{\eta_0}\) the probability measure \(\P\), given the vertex locations \(\eta_0\).
Note that, given the vertex locations, \(\P^{\eta_0}\) is a product measure with \(\operatorname{Uniform}(0,1)\) marginals and that all these events are decreasing when vertex and/or edge marks are decreased (and therefore new edges are added). Thus, applying the FKG-inequality~\cite{FKG71}, we infer
\begin{equation*}
	\P^{\eta_0}\Big(\bigcap_{k\in\N}\set{\chi(k)=0}\Big) \geq \prod_{k\in\N}\P^{\eta_0}\set{\chi(k)=0}.
\end{equation*}
To show that the product on the right-hand side is bounded away from zero, it suffices to show the equivalent statement that 
\[
\sum_{k\in\N} \P^{\eta_0}\set{\chi(k)= 1} <\infty. 
\]
which in particular implies that the probability in~\eqref{eq:noCrossing} is strictly larger than zero if
\[
\sum_{k\in\N}\mathbf{E}_0\big[\P^{\eta}\{\chi(k)=1\}\big]=\sum_{k\in\N}\P\{\chi(k)=1\}<\infty.
\]
For \(k\geq 2\), we have by symmetry,
\begin{equation*}
	\begin{aligned}
		\P\set{\chi(k)= 1}	 
		& = \P\set[\big]{\Gamma_k^{rr}\sim\Gamma_k^{\ell\ell}}+2\P\set[\big]{\Gamma_k^{\ell\ell}\sim\Gamma_k^r}\leq 3 \P\set[\big]{\Gamma_k^{\ell\ell}\sim\Gamma_k^r}.
	\end{aligned}
\end{equation*} 
The following lemma shows that for profile functions satisfying \(\rho(0+)<1\) the probability on the right hand side is bounded by the term bounded in Assumption~\eqref{eqFinite}, which immediately implies Proposition~\ref{thmFiniteComponents}.

\begin{lemma} Assume that \(\rho\) satisfies \(\rho(0+)<1\). Then for all \(\beta>0\), there exist constants \(c>0\) and \(K\in\N\) such that for all \(k\geq K\), we have
	\[
	\P\set{\Gamma_k^{\ell\ell}\sim\Gamma_k^r}\leq c 2^{2k}\int_{2^{-(1+\mu)k}}^1 \d s \int_{2^{-(1+\mu)k}}^1 \d t \ \rho(\beta^{-1}g(s,t)2^k).
	\]
\end{lemma}		

\begin{proof}
	We begin by modifying the definition of \(\mu\)-regularity (cf.\ Section~\ref{sec:LongConnect}) since we are now interested in upper bounds on connection probabilities. Throughout this proof, we say a set \(\Gamma_k^o\), \(o\in\{\ell,\ell\ell, r,rr\}\), is \emph{\(\mu\)-regular} if, for all \(i\in \set{1,\dots, \lceil 2^{k(1-\mu)}\rceil}\),
	\begin{enumerate}
		\item[(i)] \(\sum\limits_{T:(X,T)\in\Gamma_k^o} \mathbbm{1}_{\set{T\leq \lceil 2^{-(1+\mu)k}\rceil}}=0\), 
		\item[(ii)] \(\sum\limits_{T:(X,T)\in\Gamma_k^o}\mathbbm{1}_{\set{T\leq \nicefrac{i}{\lceil 2^{(1-\mu)k}\rceil}}}\leq \frac{i 2^{k+1}}{\lceil 2^{(1-\mu)k}\rceil}.\)
	\end{enumerate}
	Note, that the event of Assumption~(i) occurs with extremely high probability as the complementary event has probability  \(O(2^{-k\mu})\) for large \(k\). For Assumption~(ii), we use Chernoff's bounds to deduce 
	\[
	\P\Big\{\sum_{T:(X,T)\in\Gamma_k^o}\mathbbm{1}_{\set{T\leq \nicefrac{i}{\lceil 2^{(1-\mu)k}\rceil}}}> \frac{i 2^{k+1}}{\lceil 2^{(1-\mu)k}\rceil}\Big\}\leq 2^{(1-\mu)k}\exp(-c 2^{\mu k}).
	\]	
	Applying a straightforward union bound, it follows that the event
	\[\bigcap_{o\in\set{\ell,\ell\ell, r, rr}} \set{\Gamma_k^o \text{ is }\mu \text{-regular}}\]
	occurs with probability at least \(1-\varepsilon\) for any \(\varepsilon\in(0,1)\), if $k$ is sufficiently large. 
	{
		Furthermore, distances between end vertices of a left and right box cannot be too small due to the renewal structure. More precisely, by the law of large numbers for renewal processes, we have
		\begin{equation}\label{eq:LLNforRP}
			\frac{\sharp(\eta\cap[-n,n])}{2n}\longrightarrow \lambda, 
		\end{equation}
		\(\mathbf{P}\)-almost surely, where \(\lambda>0\) denotes the intensity of the process. Thus, we obtain for some sufficiently small (\(\lambda\) dependent) constant \(a>0\) and correspondingly chosen \(\epsilon>0\),
		\[
		\mathbf{P}_0\{|X_{-2^{k+1}}-X_{2^k}|<a \, 2^k\}=\mathbf{P}_0\big\{\sharp(\eta\cap[-n,n])>3(1+\epsilon)2^k\big\}\longrightarrow 0,
		\]
		as \(k\to\infty\).
	}   
	We thus deduce that the event
	\[
	\cE_k:=\Big(\bigcap_{o\in\set{\ell,\ell\ell, r, rr}} \set{\Gamma_k^o \text{ is }\mu \text{-regular}}\Big) \cap \set{|X_{-2^{k+1}}- X_{2^k}|>a 2^k}
	\]
	occurs with probability at least \(1-\varepsilon\), for large enough \(k\).
	
	We now argue as in the proof of Lemma~\ref{lemLongConnection}. We adapt the notion of properly spaced configurations and call such a point configuration $\omega=(x_i\colon i\in\Z)$ \emph{properly spaced} if it satisfies \(|x_{-2^{k+1}}-x_{2^k}|>a 2^k\); again we suppress the dependence of \(k\) in the notation. The property of $\mu$-regularity is measurable with respect to vertex marks only and thus independent of the vertex locations and edge marks. Hence, given a fixed configuration of properly spaced vertex locations $\omega$, we obtain
	\begin{equation}\label{eqNoEdge}
		\begin{aligned} 
			{\E}^{\eta_0=\omega} 
			& 
			\big[\1_{\{\Gamma_k^{\ell\ell}\not\sim\Gamma_k^r\}}\1_{\cE_k}\big] 
			\\ & 
			\geq {\E}^{\eta_0=\omega} \Bigg[\1_{\cE_k}\prod_{\substack{(x,T)\in\Gamma_k^{\ell\ell} \\ (y,S)\in\Gamma_k^r}} (1-\rho\big(\beta^{-1} g(S,T)a 2^k\big)\Bigg] 
			\\ & 
			\geq {\E}^{\eta_0=\omega}\Bigg[\1_{\cE_k} \exp\Big(-c \sum_{\substack{(x,T)\in\Gamma_k^{\ell\ell} \\ (y,S)\in\Gamma_k^r}}\rho(\beta^{-1}g(S,T) a 2^k)\Big)\Bigg] 
			\\	& 
			= {\E}^{\eta_0=\omega}\Bigg[\1_{\cE_k} \exp\Big(-c  \int_0^1 2^k F_{\Gamma_{k}^{\ell\ell}}(\d s)\int_0^1 2^k F_{\Gamma_k^{r}}(\d t)\ \rho(\beta^{-1}g(s,t) a 2^k)\Big)\Bigg],  
		\end{aligned}
	\end{equation}
	for some constant \(c>0\), where the second to last inequality follows from the fact that \(\rho(x)\leq \rho(0+)<1\) for all \(x>0\) by assumption. Here, \(F_{\Gamma_k^{\ell\ell}}\) denotes the empirical distribution function of the vertex marks in \(\Gamma_k^{\ell\ell}\). By \(\mu\)-regularity of the marks, we have by a similarly argument as used to derive \eqref{eqEmpDF} that
	\begin{equation*}
		\begin{aligned}
			2^k F_{\Gamma_{k}^{\ell\ell}}(t)  
			& 
			\leq 2^k \sum_{j=1}^{\lceil 2^{k(1+\mu)}\rceil} \frac{2j}{\lceil 2^{k(1+\mu)}\rceil} \mathbbm{1}_{\set{j-1<t\lceil 2^{k(1+\mu)}\rceil\leq j}} 
			\leq \frac{\big\lceil t \lceil 2^{k(1+\mu)}\rceil \big\rceil}{\lceil 2^{k(1+\mu)}\rceil} 2^k 
			\\ &
			\leq c' 2^k(t+2^{-k(1+\mu)})
		\end{aligned}
	\end{equation*}
	for some \(c'\geq 2\) uniformly for all $\mu$-regular vertex mark configurations. Plugging this into \eqref{eqNoEdge}, we infer
	\begin{equation*}
		{\P}^{\eta_0=\omega}(\{\Gamma_k^{\ell\ell}\not\sim\Gamma_k^r\}\cap \cE_k) \geq  \exp\Big(-c 2^{2k} \int_{2^{-k(1+\mu)}}^{1} \ d s \int_{2^{-k(1+\mu)}}^1 \d t \ \rho\big(\beta^{-1}g(s,t)a 2^k\big)\Big){\P}^{\eta_0=\omega}(\cE_k).
	\end{equation*}
	Since this bound is uniform in the properly spaced configuration $\omega$ for large enough \(k\), we conclude
	\begin{equation*}
		\begin{aligned}
			\P\set{\Gamma_k^{\ell\ell}\sim\Gamma_k^r} & \leq 1 -\exp\Big(-c 2^{2k} \int_{2^{-k(1+\mu)}}^{1} \ d s \int_{2^{-k(1+\mu)}}^1 \d t \ \rho\big(\beta^{-1}g(s,t)a 2^k\big)\Big) \\
			& \asymp 2^{2k} \int_{2^{-k(1+\mu)}}^{1} \ d s \int_{2^{-k(1+\mu)}}^1 \d t \ \rho\big(\beta^{-1}g(s,t)a 2^k\big),
		\end{aligned}
	\end{equation*}
	as desired.
\end{proof}
\begin{remark} \label{rem:ProofFiniteComp} 
	Besides the ergodicity induced by the i.i.d.\ nature of the location differences, we only used that \(\eta\) is a renewal process when we applied the law of large numbers in~\eqref{eq:LLNforRP} in order to justify the crossing edge argument. However, if one replaces the renewal process by any stationary and ergodic simple point process, one can easily replace the law of large numbers by the mean ergodic theorem for point processes~\cite[Theorem~8.14]{LastPenrose2017} to obtain the same result. This justifies Remark~\ref{rem:finiteComp}.
\end{remark}

It remains to prove Corollary~\ref{CorolFiniteCompPois}, i.e., the statement that the assumption \(\rho(0+)<1\) can be dropped when the vertex locations are given by a standard Poisson process.

\begin{proof}[Proof of Corollary \ref{CorolFiniteCompPois}] 
	Let \(\eta\) be a Poisson point process of intensity \(\lambda>0\). In this case, \(\beta\) can be seen as a scaling parameter of the Euclidean distance between the vertices and therefore varying \(\beta\) is equivalent to varying the intensity of the Poisson process. To see this, one can perform a linear coordinate transform on the underlying space and applying the mapping theorem for Poisson processes~\cite[Theorem~5.1]{LastPenrose2017}. 
	
	We now fix an arbitrarily \(\beta>0\) and show that no infinite component exists in \(\cG_\beta\) constructed on the Poisson process \(\eta\), or rather on its Palm version \(\eta_0\). By Poisson thinning \cite[Corollary~5.9]{LastPenrose2017}, we can interpret \(\cG_\beta\) as the graph resulting from i.i.d.\ Bernoulli site percolation of the graph \(\cG_{\nicefrac{\beta}{p}}\) for some arbitrary \(p<1\). That is, each site and all its adjacent edges are independently removed from \(\cG_{\nicefrac{\beta}{p}}\) with probability \(1-p\). We perform Bernoulli bond percolation on the graph \(\cG_{\nicefrac{\beta}{p}}\) with retention parameter \(p'\in(p,1)\) i.e., each edge is independently removed with probability \(1-p'\). By construction, this coincides with constructing the graph \(\cG_{\nicefrac{\beta}{p}}\) with the profile function \(\rho\) replaced by \(p'\rho\). Hence, we are working with the graph \(\cG_{\nicefrac{\beta}{p},\rho,g}(\xi_0)\), its bond percolated version \(\cG_{\nicefrac{\beta}{p},p'\rho,g}(\xi_0)\) and its site percolated version \(\cG_{\beta,\rho,g}(\xi_0)\). Since site percolation removes at least as many edges from the graph as bond percolation, see e.g.\ \cite{GrSt98}, and \(p'>p\), we have
	\[\P\set{\0\leftrightarrow\infty \text{ in } \cG_{\beta,\rho,g}(\xi_0)}\leq \P\set{\0\leftrightarrow\infty\text{ in }\cG_{\nicefrac{\beta}{p},p'\rho,g}(\xi_0)}.\]
	Note that assumption \eqref{eqFinite} is still satisfied and that we additionally have \(p'\rho(0+)<1\). Hence, the right hand side equals zero by Proposition~\ref{thmFiniteComponents}, finalising the proof. 
\end{proof}

\bigskip

\noindent\textbf{\large Acknowledgement.} We gratefully received support by Deutsche Forschungsgemeinschaft (DFG, German Research Foundation) – grant no. 443916008 (SPP 2265) and by the Leibniz Association within the Leibniz Junior Research Group on \textit{Probabilistic Methods for Dynamic Communication Networks} as part of the Leibniz Competition.

\end{spacing}
\section*{References} 
\renewcommand*{\bibfont}{\footnotesize}
\printbibliography[heading = none]

\end{document}